\documentclass[11pt,a4paper]{smfart} 
\usepackage{mathsmf}
\usepackage{smfthm}

\usepackage{amsfonts}
\usepackage{amssymb}
\usepackage{upgreek}
\usepackage{graphicx}
\usepackage{color}
\usepackage[utf8]{inputenc}
\usepackage[T1]{fontenc}
\usepackage[french]{babel}

\providecommand{\Ker}{\mathop\mathrm{Ker}\nolimits}
\providecommand{\Ima}{\mathop\mathrm{Im}\nolimits}
\providecommand{\Id}{\mathop\mathrm{Id}\nolimits}
\providecommand{\Ric}{\mathop\mathrm{Ric}\nolimits}
\providecommand{\Vol}{\mathop\mathrm{Vol}\nolimits}
\providecommand{\Hess}{\mathop\mathrm{Hess}\nolimits}
\providecommand{\Tr}{\mathop\mathrm{Tr}\nolimits}
\providecommand{\diam}{\mathop\mathrm{diam}\nolimits}
\providecommand{\injrad}{\mathop\mathrm{injrad}\nolimits}
\providecommand{\dimension}{\mathop\mathrm{dim}\nolimits}
\providecommand{\codimension}{\mathop\mathrm{codim}\nolimits}
\providecommand{\Det}{\mathop\mathrm{Det}\nolimits}
\providecommand{\SO}{\mathop\mathrm{SO}\nolimits}
\providecommand{\de}{\mathop\mathrm{ d }\nolimits}

\providecommand{\N}{\mathbb N}
\providecommand{\Z}{\mathbb Z}
\providecommand{\R}{\mathbb R}

\title{Petites valeurs propres des fibrés principaux en tores}
\author{Pierre Jammes}
\address{Université Côte d’Azur, CNRS, LJAD, France}

\email{pjammes@unice.fr}

\begin{document}
\begin{abstract}
Soit $M^n$ un fibré principal en tores $T^k$ sur une variété compacte
$N$. On étudie le spectre du laplacien de Hodge-de~Rham lors d'effondrements 
de $M$ sur $N$ tels que la courbure 
sectionnelle et le diamètre de $M$ vérifient $|K(M)|\leq a$ et $\diam(M)<d$.
On montre d'une part que pour tout $k$, il existe des effondrements
pour lesquels la première valeur propre du laplacien agissant sur
les formes différentielles de degré $1$ et $2$ est de l'ordre de 
$\textrm{inj}(M)^{2k}$, et d'autre part que la première valeur propre 
non nulle du laplacien agissant sur les $1$-formes est minorée
par $c(n,a,d,N)\cdot\Vol(M)^2$ et $c\cdot \textrm{inj}(M)^{2k}$
quand $M^n$ s'effondre sur $N$.

\end{abstract}
\begin{altabstract}
Let $M^n$ be a compact manifold of dimension $n$ with free $T^k$-action. 
We consider collapsings of $M$ on $N=M/T^k$ such that the sectional 
curvature and diameter of $M$ satisfy $|K(M)|\leq a$ and $\diam(M)<d$, 
and give examples of collapsings for all $k$ such that the 
first non-zero eigenvalue of Laplacian acting on $1$-forms and $2$-forms 
of $M$ are bounded above by $c(M)\cdot\textrm{inj}(M)^{2k}$.
Moreover we prove that the first non-zero eigenvalue of Laplacian 
acting on $1$-forms of all principal $T^k$-bundle $M$ over $N$ is bounded 
below by $c(n,a,d,N)\cdot\Vol(M)^2$ and $c\cdot \textrm{inj}(M)^{2k}$
when $M$ collapses on $N$.
\end{altabstract}


\maketitle

\section{Introduction}
Soit $(M,g)$ une variété riemannienne compacte connexe orientable
de dimension $n$. 
On considère l'opérateur $\Delta=\de\updelta+\updelta\de$ agissant sur l'espace 
$\Omega^p(M)$ des $p$-formes différentielles sur $M$. Le spectre de cet 
op\'erateur forme un ensemble discret de nombres positifs ou nuls qu'on
notera
$$0=\lambda_{p,0}(M,g)<\lambda_{p,1}(M,g)\leq\lambda_{p,2}(M,g)\leq\dots,$$
où la multiplicité de $\lambda_{p,0}(M,g)$ est le $p$-ième nombre de
Betti de $M$, les autres valeurs propres étant répétées s'il y a 
multiplicité.

On sait qu'à diamètre borné et courbure de Ricci minorée, la première 
valeur propre non nulle du laplacien agissant sur les fonctions est
uniformément minorée (\cite{gr80}, \cite{ly80}). B.~Colbois et G.~Courtois 
ont montré dans \cite{cc90} qu'il n'en va pas de même pour les formes 
différentielles de degré quelconque et qu'à diamètre et courbure 
sectionnelle bornés, on peut avoir des valeurs propres arbitrairement 
petites. De plus, ils montrent que ce phénomène de petites valeurs propres
est lié au fait que la variété s'effondre, c'est-à-dire que son volume 
---~ou de manière équivalente son rayon d'injectivité~---
tend vers zéro, ce qui conduit à considérer le problème suivant :
\begin{enonce}[remark]{Question}\label{intro:qu1}
Peut-on minorer $\lambda_{p,1}$ en fonction de bornes sur le diamètre,
la courbure et le rayon d'injectivité, notamment lorsque la variété 
s'effondre ?
\end{enonce}

T.~Mantuano donne dans \cite{ma08} une minoration très générale 
du spectre du laplacien de Hodge-de~Rham, mais sa dépendance par
rapport au rayon d'injectivité n'est pas explicite. S.~Chanillo et F.~Trèves 
avaient auparavant donné dans \cite{ct97} une minoration qui permettait
d'expliciter le rôle du rayon d'injectivité, mais leur résultat
est compromis par une erreur décelée par T.~Mantuano.
Dans \cite{cc00}, B.~Colbois et G.~Courtois étudient le cas particulier
des fibrés en cercles s'effondrant sur leur base et montrent ---~entre 
autres choses~--- que les petites valeurs propres sont alors de l'ordre 
du rayon d'injectivité au carré, qui est alors du même ordre que le volume
au carré. On peut donc reformuler la question \ref{intro:qu1} ainsi :
\begin{enonce}[remark]{Question}\label{intro:qu2}
Peut-on obtenir une minoration de $\lambda_{p,1}(M,g)$ asymptotiquement
de l'ordre de $\injrad(M,g)^2$ ou $\Vol(M,g)^2$ quand la variété s'effondre ?
\end{enonce}

 En dehors du cas des fibrés en cercles, cette question est presque 
totalement ouverte. Pour les situations
de limites adiabatiques sur certains feuilletages (c'est-à-dire lorsque la 
métrique s'écrit sous la forme $g=g_H\oplus\varepsilon^2g_V$ où
$g_V$ une métrique sur les espaces tangents $V$ aux feuilles, $g_H$
une métrique sur une distribution complémentaire à $V$), on sait
(\cite{fo95}, \cite{alk00}) qu'on peut calculer pour tout entier $k$ le nombre 
de petites valeurs propres de l'ordre de $\varepsilon^{2k}$ à l'aide
d'une suite spectrale, mais ces effondrements ne sont pas nécessairement
à courbure bornée. Dans le cas où la feuille est de dimension~1,
on peut cependant obtenir des effondrements à courbure bornée produisant
des petites valeurs propres (cf.~\cite{ja10}). On peut trouver une 
présentation plus détaillée des questions~\ref{intro:qu1} et~\ref{intro:qu2} 
dans~\cite{ja05} ou~\cite{jahdr}.

Le but du présent travail est de donner des éléments de réponse à
la question \ref{intro:qu2} en étudiant le spectre 
des fibrés principaux en tores s'effondrant sur leur
base. Nous allons montrer deux résultats. Le premier apporte une réponse
négative à la question \ref{intro:qu2} en ce qui concerne le rayon
d'injectivité. En effet on peut construire pour tout $k\in\N^*$ des 
effondrements pour lesquels certaines valeurs propres décroissent
au moins aussi vite que $\injrad(M,g)^{2k}$. Pour cela on utilisera
le fait que, si la fibre est de dimension supérieure ou égale à 2, on
peut l'effondrer non par homothétie mais en privilégiant une
certaine direction qu'on choisira pour ses propriétés diophantiennes.
Ces effondrements permettent aussi de construire des contre-exemples
à certaines des estimées de \cite{ct97}:
\begin{theo}\label{intro:th5}
Pour tout entier $k\geq1$ et pour toute variété $(N,h)$
telle que $b_2(N)\geq k$, il existe un fibré principal $M$ en tore $T^k$
sur $N$, une famille de métrique $(g_\varepsilon)_{\varepsilon\in]0,1]}$
sur $M$, et des réels strictement positifs $C_i(k,(N,h))$, $i=1,2,3$
et $\varepsilon_0(k,(N,h))$
tels que la courbure et le diamètre de $(M,g_\varepsilon)$ soient
uniformément bornés par rapport à $\varepsilon$,
$\Vol(M,g_\varepsilon)=\varepsilon$ pour tout $\varepsilon$, et
\begin{equation}\label{intro:eq3}
\varepsilon^2\leq C_1\cdot\lambda_{p,1}(M,g_\varepsilon)\leq 
C_2\cdot\injrad^{2k}(M,g_\varepsilon)\leq C_3\cdot\varepsilon^2
\end{equation}
pour $p=1$ et 2, et pour tout $\varepsilon<\varepsilon_0$.

 De plus, si $b_1(N)>b_2(M)$, on a aussi pour $p=2$ et 3
\begin{equation}\label{intro:eq4}
\varepsilon^2\leq C_1\cdot\lambda_{p,b_1(N)-b_2(M)}(M,g_\varepsilon)\leq 
C_2\cdot\injrad^{2k}(M,g_\varepsilon)\leq C_3\cdot\varepsilon^2.
\end{equation}
\end{theo}
\begin{rema}
Le théorème \ref{intro:th5} montre qu'on ne peut pas minorer la première
valeur propre du laplacien par une puissance de rayon d'injectivité
indépendante de la dimension. Plus précisément, comme on peut construire
des variétés de dimension 3 dont le 2\up{e} nombre de Betti est
arbitrairement grand, il existe pour tout $n\geq4$ des variétés de
dimension $n$ vérifiant (\ref{intro:eq3}) avec $k=n-3$. On ne peut donc
pas espérer obtenir de minoration générale par le rayon d'injectivité
avec un exposant inférieur à $2n-6$ si $n\geq4$.
\end{rema}

\begin{rema}\label{intro:rem}
Les effondrements du théorème~\ref{intro:th5}
fournissent des contre-exemples à certains des résultats  de~\cite{ct97},
en particulier les théorèmes~1.1 et~1.2 (voir exemple~\ref{inj:ex2}).
\end{rema}

\begin{rema}
On verra dans la section \ref{inj} qu'il existe effectivement
des fibrés vérifiant la condition $b_1(N)>b_2(M)$.
\end{rema}

En petite dimension, on peut en déduire un résultat pour les formes
de tout degré :
\begin{coro}\label{intro:cor}
Pour $n=5$ et $7$, il existe un fibré $M$ en tores $T^2$ sur une variété
de dimension $n-2$, une suite de métrique $(g_\varepsilon)$ sur $M$ et 
une constante $C>0$ telles que 
la suite $(g_\varepsilon)$ effondre le fibré $M$ sur sa base à courbure
et diamètre bornés quand $\varepsilon$ tend vers zéro, et que
\begin{equation}\label{intro:eq5}
\lambda_{p,1}(M,g_\varepsilon)\leq C\cdot\injrad^{4}(M,g_\varepsilon)
\end{equation}
pour $1\leq p\leq n-1$.
\end{coro}

Le second théorème donne, dans le cas des $1$-formes, 
des minorations du spectre en fonction du volume et du rayon d'injectivité :
\begin{theo}\label{intro:th4}
Soit deux réels $a$ et $d$ strictement positifs, un entier $n\geq3$ et
$(N,h)$ une variété riemannienne de dimension strictement inférieure 
à $n$. Il existe des constantes $\varepsilon_0(n,a,d,(N,h))>0$, 
$C(n,a,d,(N,h))>0$ et $C'(n,a,d,(N,h))>0$ telles que si $(M,g)$ 
est une variété riemannienne 
de dimension $n$ vérifiant $\diam(M,g)\leq d$, $|K(M,g)|\leq a$ et 
si $\pi:(M,g)\rightarrow(N,h)$ est une fibration principale de fibre 
$T^k$ qui soit une $\varepsilon$-approximation de Hausdorff avec 
$\varepsilon<\varepsilon_0$, alors
\begin{equation}\label{intro:eq2}
\lambda_{1,1}(M,g)\geq C\cdot\Vol^2(M,g)\geq C'\cdot\injrad^{2k}(M,g).
\end{equation}
\end{theo}

\begin{rema}
L'exposant du volume dans la minoration (\ref{intro:eq2}) est optimal. 
On peut par exemple
considérer le cas de produit riemannien d'un fibré en cercles qui s'effondre
par un tore dont la métrique est fixée mais dont le diamètre est
suffisamment petit pour que le fibré produit soit proche de la base pour
la distance de Gromov-Hausdorff, les résultats de \cite{cc00} assurent
que dans ce cas, $\lambda_{1,1}(M)$ ne peut pas décroître plus lentement 
que $\Vol^2(M)$. L'inégalité (\ref{intro:eq3}) assure par ailleurs
que l'exposant du rayon d'injectivité dans (\ref{intro:eq2}) est lui
aussi optimal.
\end{rema}

\begin{rema}
On sait déjà que si la fibre est de dimension supérieure ou égale à deux,
on ne peut pas obtenir pour les fibrés en tores de majoration
de $\lambda_{1,1}(M)$ par le volume au carré comme c'est le cas pour 
les fibrés en cercles. En effet,
on peut construire des exemples d'effondrements de fibrés principaux en tores
pour lesquels la première valeur propre du laplacien ne tend
pas vers zéro (cf. \cite{ja03}, théorème 1.13 et paragraphe 4.2).
\end{rema}

Le théorème \ref{intro:th4} soulève les deux questions suivantes qui
restent ouvertes:
\begin{enonce}[remark]{Question} Peut-on généraliser ces résultats 
aux $p$-formes différentielles, pour tout $p$ ?
\end{enonce}
\begin{enonce}[remark]{Question}
La minoration du spectre par le volume au carré de la variété se
généralise-t-elle à d'autres familles de variétés ?
\end{enonce}

Pour démontrer les théorèmes \ref{intro:th5} et \ref{intro:th4}, nous 
commencerons par trois sections préliminaires sur la topologie, la
 géométrie et le spectre des fibrés principaux en tores.

 Dans la section \ref{topo}, nous définirons un invariant 
topologique généralisant la classe d'Euler des fibrés en cercles et qui 
pourra être utilisé pour contrôler le spectre.

 La section \ref{geom} sera consacrée à l'étude de la géométrie des
fibrés principaux en tores. Le résultat principal est que pour
démontrer le théorème \ref{intro:th4} on
peut se ramener à une situation géométrique simple.
En particulier, on montrera qu'une métrique de courbure et diamètre
bornés sur le fibré est proche d'une métrique invariante pour laquelle
les fibres sont totalement géodésiques (théorème~\ref{geom:th}).

 Dans la section \ref{inv}, nous étudierons comment, dans le cas d'un fibré 
principal en tore $T^k$ muni d'une métrique invariante, on peut se
ramener à l'étude des petites valeurs propres du laplacien à celle
du spectre du laplacien resteint aux formes différentielles
invariantes par l'action de $T^k$ (théorème~\ref{inv:th2}). Dans
le cas des fibrés en cercles, ce théorème  améliore
sensiblement un résultat semblable obtenu par B.~Colbois et G.~Courtois
(\cite{cc00}) qui nécessitait des hypothèses de courbure et diamètre
bornés.

 Les deux dernières sections serons consacrées aux démonstrations des
théorèmes \ref{intro:th5} et \ref{intro:th4}.
 Dans la section \ref{inj}, nous construirons les exemples d'effondrements
annoncés dans le théorème \ref{intro:th5}.

Enfin, dans la section \ref{adap}, nous démontrerons l'inégalité
(\ref{intro:eq2}) du théorème~\ref{intro:th4} en utilisant les résultats 
des sections précédentes.

 Ces travaux ont été en grande partie réalisés à l'Institut de mathématiques 
de Neuchâtel et au Laboratoire de mathématiques d'Avignon, qu'ils en
soient remerciés. Ma gratitude va aussi à Amandine Lacuès pour sa 
relecture attentive du texte, ainsi qu'au rapporteur anonyme dont les
remarques ont permis d'améliorer le texte de de simplifier certaines
démonstrations.

\section{Topologie des fibrés principaux en tores}\label{topo}
Nous allons dans cette partie nous attacher à décrire la topologie
des fibrés principaux en tores, et en particulier à construire
un invariant différentiel qui permettra, comme la classe d'Euler
dans la cas des fibrés en cercles, d'étudier le comportement du spectre
du laplacien lors d'un effondrement.

Soit $M$ un fibré principal en tores $T^k$ sur une base $N$. Le
fibré $M$ peut se décomposer en la somme de Whitney de $k$
fibrés en cercles, de classes d'Euler respectives $e_1,\ldots,e_k$, 
mais le $k$-uplet $(e_1,\ldots,e_k)\in H^2(N)^k$ dépend du choix de 
cette décomposition.

Dans le cas d'un fibré en cercles de classe d'Euler $[e]$, on a la 
propriété suivante (\cite{bt82}, p.72)~: si $\omega$ est une $1$-forme
de connexion du fibré, alors $\de\omega$ est
une $2$-forme horizontale qui dépend du choix de la connexion sur le fibré, 
mais qui est, au signe près, le relevé d'un élément de $[e]$. Dans le
cas d'un fibré en tores, on va construire un invariant qui généralise
cette propriété. 

 Dans le cas d'un fibré principal quelconque, une forme de connexion 
$\gamma$ est une $1$-forme verticale invariante à valeur dans l'algèbre 
de Lie $\mathcal G$ de la fibre (voir par exemple~\cite{ra05}, chapitre~5). 
Si on se donne un élément $\omega$ de 
$\mathcal G^*$ et qu'on applique $\omega$ à l'image de la forme de 
connexion, on obtient une $1$-forme verticale $\omega\circ\gamma$ à 
valeur réelle. Si $\Ker\omega$ est tangent à un sous-tore $T^{k-1}
\subset T^k$, alors $\omega\circ\gamma$ est le relevé d'une forme de connexion
sur $M/T^{k-1}$, et donc $\de(\omega\circ\gamma)$ est le relevé de la
classe d'Euler $e\in H^2(N)$ de ce fibré. Par linéarité, 
on peut étendre cette construction en une application $\omega\to e(\omega)$.

\begin{definition}
On appelle classe d'Euler du fibré principal $T^k\hookrightarrow M
\stackrel\pi\to N$
l'application $e: \mathcal G^*\to H^2(N,\R)$ déterminée par 
$\de(\omega\circ\gamma)=\pi^*(e(\omega))$ où $\gamma$ est une forme de
connexion du fibré et $\mathcal G$ l'algèbre de Lie de $T^k$.
\end{definition}
\begin{rema} 
Si $k=1$ et si $\omega$ est la forme volume du cercle de longueur $1$,
alors $e(\omega)$ est la classe d'Euler du fibré.
\end{rema}

 On a défini ainsi une classe d'Euler en cohomologie de~Rham. On peut
aussi la définir en cohomologie de \v{C}ech en partant du fait 
que les $T^k$-fibrés principaux ont classifiés par les éléments
de $H^1(N,T^k)$. Comme la suite exacte courte $1\to\Z^k\to\R^k\to T^k\to1$
induit une suite exacte en cohomologie de \v{C}ech
\begin{equation}
H^1(N,\Z^k)\to H^1(N,\R^k)\to H^1(N,T^k)\to H^2(N,\Z^k)\to,
\end{equation}
on peut définir la classe d'Euler comme l'image dans $H^2(N,\Z^k)$ de 
l'élément de $H^1(N,T^k)$ caractérisant le fibré. L'équivalence entre 
les deux définitions s'établit en se ramenant par projections au cas des
fibrés en cercles.

\begin{exem}
Le théorème 1.13 de \cite{ja03} montre ---~entre autres choses~--- 
qu'un fibré 
principal en tores $T^k$ non trivial dont la base est un tore $T^2$ peut 
s'écrire comme produit d'une nilvariété de dimension 3 et d'un tore 
$T^{k-1}$. Comme $H^2(T^2)$ est 
de dimension $1$, le noyau de $e$ est de dimension $k-1$, ce qui signifie
qu'on peut décomposer le fibré en une somme de Whitney de $k$ fibrés
en cercles dont $k-1$ sont triviaux. On retrouve donc le fait que
le fibré peut s'écrire comme le produit d'un fibré en cercles sur
$T^2$ et d'un tore de dimension $k-1$.
\end{exem}

\section{Géométrie des fibrés principaux en tores}\label{geom}
\subsection{Métriques adaptées}
Nous allons ici montrer qu'on peut, dans le but d'obtenir le
théorème~\ref{intro:th4}, se ramener à une situation géométrique 
plus simple.
Cette situation est une généralisation de la notion de métrique 
adaptée définie dans le cas des fibrés en cercles par B.~Colbois
et G.~Courtois (\cite{cc00}) :
\begin{definition}\label{geom:adap}
On dit que le couple de métriques $(g,h)$ définies sur $M$ et $N$ 
respectivement est adapté à la fibration principale
$T^k\hookrightarrow M^n\stackrel{\pi}{\rightarrow} N$ si:
\begin{enumerate}\label{geom:df}
\item La fibration $\pi:(M,g)\rightarrow (N,h)$ est une submersion 
riemannienne ;
\item L'action de $T^k$ sur $M$ est isométrique ;
\item Les fibres sont totalement géodésiques ;
\item Toute $1$-forme verticale $\omega$ induite par un élément de 
$\mathcal G^*$ vérifie $\de\omega=\pi^*(e(\omega))$, où $e$ est
la classe d'Euler du fibré.
\end{enumerate}
\end{definition}

On veut montrer qu'une métrique de courbure bornée
sur un fibré principal en tores est proche d'une métrique 
adaptée:
\begin{theo}\label{geom:th}
 Soient $a$ et $d$ deux réels strictement positifs, et
$T^k\hookrightarrow (M^n,g)\stackrel{\pi}{\rightarrow} (N,h)$
un fibré principal en tores. Il existe 
des constantes $\varepsilon_0(n,a,(N,h))>0$, $\tau(n,a,(N,h))>0$, 
$\tau'(n,a,(N,h))>0$ et $c(n,a,(N,h))>0$
telles que si $|K(N,h)|\leq a$, $|K(M,g)|\leq a$ et
si $\pi$ est une $\varepsilon$-approximation de Hausdorff avec $\varepsilon
<\varepsilon_0$, alors il 
existe des métriques $\tilde g$ et $\tilde h$  sur $M$ et $N$ respectivement 
et une fibration $\pi': (M,\tilde g)\rightarrow(N,\tilde h)$  telles que 
\begin{enumerate}
\item Le couple $(\tilde g,\tilde h)$ est adapté à la fibration $\pi'$;
\item $\displaystyle\frac 1\tau g\leq\tilde g\leq \tau g$ et
$\displaystyle\frac 1\tau h\leq\tilde h\leq \tau h$;
\item La restriction de $\tilde g$ à la fibre est telle que 
$\diam(\pi'^{-1}(x))\leq\tau'\varepsilon$, pour tout $x\in N$;
\item La courbure sectionnelle de $(M,\tilde g)$ vérifie
$|K(X,Y)|\leq c$, pour toute paire de vecteurs horizontaux orthonormés
$(X,Y)$.
\end{enumerate}
\end{theo}
On pourra alors appliquer le résultat de J.~Dodziuk selon
lequel si deux métriques sont proches, alors les spectres du laplacien
pour ces deux métriques sont proches aussi:
\begin{theo}[\cite{do82}]
Soit $g$ et $\tilde g$ deux métriques riemanniennes sur une variété compacte
$M$ de dimension $n$, et $\tau$ une constante strictement positive. Si les deux
métriques vérifient $\frac1\tau g\leq\tilde g\leq\tau g$, alors
$$\frac1{\tau^{3n-1}}\lambda_{p,k}(M,g)\leq\lambda_{p,k}(M,\tilde g)
\leq\tau^{3n-1}\lambda_{p,k}(M,g),$$
pour tous entiers $k\geq0$ et $p\in[0,n]$.
\end{theo}

\begin{rema} En vertu d'un théorème de Hermann (\cite{her60},
\cite{be08} p.~249),
le fait que les fibres soient totalement géodésiques implique 
qu'elles sont isométriques entre elles. On va voir dans la démonstration du
théorème~\ref{geom:th} que réciproquement, sur les fibrés considérés, si
la métrique est invariante et que les fibres sont isométriques alors 
elles sont totalement géodésiques.
\end{rema}
\subsection{Situation de métrique invariante}
 Nous allons dans un premier temps montrer que si on suppose qu'on 
a sur $M$ une métrique invariante, elle est proche d'une métrique
qui vérifie les points (1) à (3) de la définition 
\ref{geom:df}. Plus précisément :
\begin{prop}\label{geom:pr}
Soit $T^k\hookrightarrow (M^n,g)\stackrel{\pi}{\rightarrow} (N,h)$
un fibré principal en tores muni d'une métrique invariante $g$ tel que 
$\pi$ soit une submersion riemannienne. Pour tout $a>0$ et $d>0$,
il existe des constantes $\tau(n,a,d)>0$ et $c(n,a)>0$ telles
que si $|K(N,h)|\leq a$, $K(M,g)\geq -a$ et $\diam(M,g)\leq d$, 
alors il existe une métrique invariante $\tilde g$ sur $M$ telle 
que la fibration $\pi:(M,\tilde g)\rightarrow(N,h)$ soit une 
submersion riemannienne à fibres totalement géodésiques et
$$\frac 1\tau g\leq\tilde g\leq \tau g.$$
\end{prop}
\begin{rema} On peut noter qu'on utilise non pas une hypothèse
de courbure bornée sur $M$ mais seulement que la courbure sectionnelle
est minorée.
\end{rema}

Pour montrer la proposition \ref{geom:pr}, on utilisera les deux lemmes 
suivants. Le premier est une application directe de la formule de O'Neill:
\begin{lem}\label{geom:oneill}
Soit $a>0$ et $T^k\hookrightarrow (M^n,g) 
\stackrel{\pi}{\rightarrow} (N,h)$ un fibré
principal en tores muni d'une métrique invariante $g$ tel que $\pi$
soit une submersion riemannienne, $|K(N,h)|\leq a$,
et $K_{(M,g)}(X,Y)\geq -a$ pour tout couple $(X,Y)$ de vecteurs 
horizontaux orthonormés.
Alors, pour toute $1$-forme différentielle 
verticale $\omega$, on a :
\begin{enumerate}
\item $\displaystyle|\de\omega(X,Y)|_x^2\leq\frac{8a}3|\omega|_x^2$, pour tout 
$x\in M$ et tout couple de vecteurs horizontaux orthonormés $X$ et $Y$ ;
\item $\displaystyle\|\de\omega\|_\infty\leq\frac{4an(n-1)}3\|\omega\|_\infty$. 
\end{enumerate}
\end{lem}
\begin{proof} Soit $x\in M$, $y=\pi(x)$, $\tilde X$ et $\tilde Y$ 
deux champs de $N$ orthonormés en $y$, et $X$ et $Y$ les relevés de
$\tilde X$ et $\tilde Y$ à $M$. La formule de O'Neill (\cite{ghl87}
p.~127, \cite{be08} p.~241) donne
\begin{equation}
K_N(\tilde X,\tilde Y)=K_M(X,Y)+\frac34\left|[X,Y]^V\right|^2,
\end{equation}
où $[X,Y]^V$ désigne la composante verticale de $[X,Y]$.
D'autre part on a, en utilisant le fait que $\omega$ est verticale, 
\begin{eqnarray}
\de \omega(X,Y)&=&X\cdot\omega(Y)-Y\cdot\omega(X)-\omega([X,Y])\nonumber\\
&=&-\omega([X,Y]).
\end{eqnarray}
On en déduit :
\begin{eqnarray}
|\de \omega(X,Y)|_x^2&=&|\omega([X,Y])|_x^2\leq|\omega|_x^2|[X,Y]^V|_x^2
\nonumber\\
&\leq&\frac43|\omega|_x^2(K_y(\tilde X,\tilde Y)-K_x(X,Y)).
\end{eqnarray}
Comme chacun des couples $(\tilde X,\tilde Y)$ et $(X,Y)$ est orthonormé
en $x$ et $y$, on a les majorations $|K_y(\tilde X,\tilde Y)|\leq a$ 
et $K_x(X,Y)\geq -a$, et donc
\begin{equation}
|\de \omega(X,Y)|_x^2\leq\frac{8a}3|\omega|_x^2.
\end{equation}
 Et comme l'inégalité précédente est vraie quel que soit le choix
de $(\tilde X,\tilde Y)$, il en découle finalement
\begin{equation}
|\de\omega|_x^2\leq\frac{4an(n-1)}3|\omega|_x^2
\leq\frac{4an(n-1)}3\|\omega\|^2,
\end{equation}
ce qui achève la démonstration.
\end{proof}

Le second lemme montre que dans le cas d'un fibré en cercles, à courbure 
bornée, la longueur des fibres varie peu d'une fibre à l'autre. 
\begin{lem}\label{geom:lem2}
Soit $S^1\hookrightarrow (M^n,g)\stackrel{\pi}{\rightarrow} (N,h)$
un fibré principal en cercles sur $N$, tel que $g$ soit invariante
et $\pi$ soit une submersion riemannienne. Pour tout $a>0$ et $d>0$, 
il existe $\tau(n,a,d)>0$ tel que si $|K(N,h)|\leq a$, $K(M,g)\geq -a$ 
et $\diam(M,g)\leq d$, alors pour tout $x,y\in N$, on a
$$\frac 1\tau l_y\leq l_x\leq \tau l_y,$$
où $l_x$ et $l_y$ désignent les longueurs des fibres au dessus de 
$\pi^{-1}(x)$ et $\pi^{-1}(y)$ respectivement.
\end{lem}
\begin{proof}
 On choisit sur le fibré $M$ une $1$-forme 
verticale $\omega$ dont l'intégrale sur chaque fibre est égale à $1$.
Soit $U$ le champ vertical induit par l'action de $S^1$ qui vérifie
$\omega(U)=1$. La norme $|U|$ de ce champ est
constante sur chaque fibre, et s'écrit $|U|=\pi^*f$, où $f$ est une
fonction sur $N$. De plus, en tout point $x$ de $N$, la norme
$f(x)$ de $U$ est égale à la longueur de la fibre au dessus de $x$.
On va montrer que $f$ est bornée en fonction de $a$ 
et $d$. Remarque : $\omega$ n'est pas la forme duale de $U$ pour
la métrique. Sa norme ponctuelle sur la fibre $\pi^{-1}(x)$
est $|\omega|=\frac1f$, et on a $U^\flat=f^2\omega$.

Soit $x\in N$, et $\tilde X$ un vecteur unitaire tangent à $N$ en $x$.
Soit $\tilde X_i$ une base orthonormée de champs de vecteurs au voisinage 
de $x$, telle que $\mathrm{D}{\tilde X_1}\tilde X_1=0$ sur ce voisinage, et 
$\tilde X_{1|x}=\tilde X$.
On relève cette base à $T^HM$ en notant $X_i=\pi^*(\tilde X_i)$ et
$X=X_1$.  Ces champs vérifient $[U,X_i]=0$.
En effet, ces crochets de Lie sont déterminés par $[U,X_i]=\mathcal L_U X_i$,
or les champs $X_i$ sont $S^1$-invariants. On notera par ailleurs $V$ 
le champ de norme~$1$ défini par $V=U/|U|$. On a alors $[V,X_i]=(\de f(X_i)/f)V$.

On va calculer la courbure sectionnelle $K(X,V)$ en fonction de
$f$ et de ses variations. On utilisera pour cela les formules de O'Neill
(cf. chapitre~9 de  \cite{be08}) qui donnent en particulier
\begin{equation}
K(X,V)=((\mathrm{D}_XT)_VV,X)-|T_VX|^2+|A_XV|^2,
\end{equation}
avec les notations suivantes :
\begin{itemize}
\item le vecteur $A_XV$ est la composante horizontale de $\mathrm{D}_XV$;
\item le vecteur $T_VX$ est la composante verticale de $\mathrm{D}_VX$;
\item le vecteur $T_VV$ est la composante horizontale de $\mathrm{D}_VV$ (et mesure
donc la courbure extrinsèque de la fibre).
\end{itemize}
On utilisera abondamment la formule de Koszul, qui caractérise la
connexion de Levi-Civita,
\begin{eqnarray}\label{geom:lc}
2\langle\mathrm{D}_{Z_1}Z_2,Z_3\rangle&=&
Z_1\cdot\langle Z_2,Z_3\rangle
+Z_2\cdot\langle Z_3,Z_1\rangle-Z_3\cdot\langle Z_1,Z_2\rangle\nonumber\\
&&+\langle[Z_1,Z_2],Z_3\rangle-\langle[Z_1,Z_3],Z_2\rangle
-\langle[Z_2,Z_3],Z_1\rangle,\hspace{3mm}
\end{eqnarray}
ainsi que l'orthogonalité de $(X_1,\cdots,X_n,U)$ qui annulera les premiers
termes du membre de droite de cette formule.
Elle donne directement 
\begin{equation}
\langle \mathrm{D}_VX,V\rangle=\langle [V,X],V\rangle=\frac{\de f(X)}f\textrm{ et }
\langle \mathrm{D}_XV,X_i\rangle=-\frac12\langle[X,X_i],V\rangle^2.
\end{equation}
On en déduit $|T_UX|^2=\left(\frac{\de f(X)}f\right)^2$ et
$|A_XV|^2=\frac14\sum_i\langle[X,X_i],V\rangle^2$.

On calcule le terme restant en partant de la formule de Leibniz :
\begin{equation}
\langle(\mathrm{D}_XT)_VV,X\rangle=\mathrm{D}_X\langle T_VV,X\rangle
-\langle T_{\mathrm{D}_XV}V,X\rangle-\langle T_V\mathrm{D}_XV,X\rangle.
\end{equation}
La formule de Koszul donne alors 
\begin{equation}
\langle T_VV,X\rangle=\langle\mathrm D_VV,X\rangle=-\langle[V,X],V\rangle
=-\frac{\de f(X)}f
\end{equation}
et donc 
\begin{equation}
\mathrm{D}_X\langle T_VV,X\rangle=\left(\frac{\de f(X)}f\right)^2-
\frac{\Hess f(X,X)}f.
\end{equation}

On obtient aussi que $\mathrm D_XV=-\frac12\sum_i\langle[X,X_i],V\rangle X_i$.
Comme $T_{X_i}V=0$ (\cite{be08}, p.~239), on en déduit d'une part que 
$T_{\mathrm{D}_XV}V=0$. Et d'autre part, comme $T_VX_i$ est vertical on a aussi
$\langle T_V\mathrm{D}_XV,X\rangle=0$. Il reste finalement :
\begin{equation}
K(X,V)=-\frac{\Hess f(X,X)}f+\frac14\sum_i\langle[X,X_i],V\rangle^2.
\end{equation}

Les derniers termes peuvent être majorés en fonction de la courbure. En
effet, on a
\begin{equation}
\langle[X,X_i],V\rangle=V^\flat([X,X_i])=f\omega([X,X_i])
=-f\de\omega(X,X_i),
\end{equation}
et donc, en vertu du lemme \ref{geom:oneill},
\begin{equation}
\langle[X,X_i],V\rangle^2\leq f^2\frac{8a}3|\omega|^2=\frac{8a}3.
\end{equation}
Par hypothèse, la courbure sectionnelle de $M$ est minorée par $-a$.
On a donc finalement :
\begin{equation}\label{geom:hess}
\frac{\Hess f(X,X)}f\leq \left(\frac{2n}3+1\right)a.
\end{equation}

Soient $x$ et $y$ deux points de $N$, et $\gamma$ une géodésique
minimisante joignant ces deux points. Notons $\mu$ la fonction
définie par 
\begin{equation}\mu(t)=\ln f\circ\gamma(t).\end{equation}
En dérivant $\mu$ par rapport à $t$, on obtient :
\begin{equation}
\mu'(t)=\frac{\de f(\gamma'(t))}{f\circ\gamma(t)}
\end{equation}
et
\begin{eqnarray}
\mu''(t)&=&-\left(\frac{\de f(\gamma'(t))}{f\circ\gamma(t)}\right)^2
+\frac{\de f(\mathrm{D}_{\gamma'(t)}\gamma'(t))+\mathrm{D}\de f(\gamma'(t),\gamma'(t))}
{f\circ\gamma(t)}\nonumber\\
&=&-\mu'(t)^2+\frac{\Hess f(\gamma'(t),\gamma'(t))}{f\circ\gamma(t)}.
\end{eqnarray}
On a donc, en vertu de la majoration (\ref{geom:hess}) :
\begin{equation}
\mu''(t)\leq\frac{\Hess f(\gamma'(t),\gamma'(t))}{f\circ\gamma(t)}
\leq\left(\frac{2n}3+1\right)a.
\end{equation}
Supposons que $x$ est un point où $f$, et donc $\mu$, atteint son 
minimum. On a alors $\mu'(0)=0$ et donc, en remarquant que $d$ majore
le diamètre de $N$,
\begin{equation}
\mu'(t)\leq\left(\frac{2n}3+1\right)at,
\end{equation}
et
\begin{equation}
\mu(t)\leq\frac{2n+3}6at^2\leq\frac{2n+3}6ad^2.
\end{equation}
Le rapport $\frac{f(y)}{f(x)}$ est donc majoré par une constante
$\tau=\exp(\frac{2n+3}6ad^2)$. Comme on a montré cette majoration 
en prenant pour $x$ un point où $f$ atteint son  minimum, elle
sera vraie \emph{a fortiori} pour un $x$ quelconque.

 Remarque : si les fibres sont isométriques, alors le champ $U$ est de
norme constante. Il est aisé de vérifier à l'aide de (\ref{geom:lc})
que $\mathrm{D}_UU$ est alors nul, c'est-à-dire que les fibres sont totalement
géodésiques.
\end{proof}

\begin{proof}[de la proposition \ref{geom:pr}]
 Le but est en fait de généraliser le lemme \ref{geom:lem2}
aux fibrés en tores pour montrer que $g$ est proche d'une métrique
pour laquelle toutes les fibres sont isométriques.

Soit $\bar U\in\mathcal G$ non nul, où $\mathcal G$ est l'algèbre de Lie de
$T^k$, et $U$ le champ vertical induit par 
$\bar U$ sur $M$. Soit $x_0\in N$. On choisit $\bar U$ de sorte que 
$|U|=1$ au dessus de $x_0$. De plus, on impose à $\bar U$ d'avoir un
coefficient directeur rationnel, c'est-à-dire que $\bar U$ est colinéaire
à un vecteur de $\Z^k\subset\mathcal G$. L'action du flot associé à
$U$ induit alors une fibration $S^1\hookrightarrow M
\stackrel{\pi}{\rightarrow}(M',g')$. On peut contrôler la courbure
de ce fibré. En effet, en utilisant la formule de O'Neill, on peut
écrire
\begin{eqnarray}
K_{M'}(\tilde X,\tilde Y)&=&K_M(X,Y)-\frac34|[X,Y]^V|^2\nonumber\\
&=&K_M(X,Y)-\frac34\frac{|\de\omega(X,Y)|^2}{|\omega|^2},
\end{eqnarray}
où $\tilde X$ et $\tilde Y$ sont deux vecteurs orthonormés de $M'$,
$X$ et $Y$ leurs relevés respectifs sur $M$, et $\omega$ la $1$-forme
induite par l'action de $T^k$ telle que 
$\omega(U)=1$. Le lemme \ref{geom:oneill} permet
de contrôler le dernier terme en fonction de la courbure de $M$, et donc
\begin{equation}
K_{M'}(\tilde X,\tilde Y)\geq -a-2a=-3a.
\end{equation}
Le lemme \ref{geom:lem2} assure alors qu'il existe une constante 
$\tau(n,a,d)$ telle que 
\begin{equation}\label{geom:eq1}
\frac1\tau\leq|U|\leq\tau.
\end{equation}
Notons $\tilde g$ la métrique invariante sur $M$ obtenue en modifiant 
$g$ dans la direction verticale de sorte que les fibres soient isométriques
à $\pi^{-1}(x_0)$, et en conservant la distribution horizontale et la
métrique horizontale associées à $g$. Pour cette nouvelle métrique, la
norme de $U$ est uniformément égale à 1.
La relation (\ref{geom:eq1}) peut s'écrire
\begin{equation}\label{geom:eq2}
\frac1{\tau^2}\tilde g(U,U)\leq g(U,U)\leq\tau^2\tilde g(U,U),
\end{equation}
Par continuité, les inégalités (\ref{geom:eq2}) s'étendent à n'importe 
quel vecteur vertical.
Comme $g$ et $\tilde g$ sont identiques sur la direction horizontale, 
on aura finalement
\begin{equation}
\frac1{\tau^2}\tilde g\leq g\leq\tau^2\tilde g.
\end{equation}

Pour conclure, remarquons que si la métrique sur $M$ est telle que les 
fibres soient isométriques, la fibration $S^1\hookrightarrow M
\stackrel{\pi}{\rightarrow}(M',g')$ induite par le champ $U$ est à fibre 
totalement géodésique. Par continuité, les fibres du fibré $T^k
\hookrightarrow M \rightarrow N$ sont elles aussi totalement géodésiques.
\end{proof}
\subsection{Cas général}

 On va maintenant démontrer le théorème \ref{geom:th}. Pour ce faire, 
on va s'inspirer d'une démonstration d'un théorème de Lott 
(\cite{lo02}, théorème 2), qui utilise les résultats de \cite{cfg92}.

 Soit $g$ une métrique sur $M$ vérifiant les hypothèses du théorème 
\ref{geom:th}. Tout d'abord, en utilisant un résultat de régularisation
d'Abresch (\cite{cfg92}, théorème 1.12), on construit une métrique
$g_1$ sur $M$ telle que $\frac1{\tau_1}g\leq g_1\leq\tau_1g$,
$|K(M,g_1)|\leq a$ et 
$\|\mathrm{D}^i R\|\leq A_i(n,a,\tau_i)$, où $\tau_1>1$ est un réel fixé, et
$\mathrm{D}$ et $R$ désignent respectivement la dérivée covariante et
le tenseur de courbure pour la métrique $g_1$.

 On applique ensuite le théorème 2.6 de \cite{cfg92}, qui assure l'existence
de constantes $\varepsilon_0(n,(N,h))$, $\kappa(n,A)$, $\kappa'(n,A,(N,h))$
et $\kappa_i(n,A,(N,h))$ et d'une fibration $\pi':(M,g_1)\rightarrow
(N,h)$ tel que si $\pi$ est une $\varepsilon$-approximation de 
Hausdorff avec $\varepsilon<\varepsilon_0(n,(N,h))$, alors :

\begin{itemize}
\item pour tout $x\in N$, le diamètre de $\pi'^{-1}(x)$ pour 
la métrique $g'$ est inférieur à $\kappa\varepsilon$;
\item la seconde forme fondamentale de la fibre vérifie 
$\|I\!I_{\pi^{-1}(x)}\|_\infty\leq\kappa'$ pour tout $x\in N$;
\item la submersion $\pi'$ est $\kappa_i$-régulière, c'est-à-dire que
$\|\mathrm{D}^i\pi'\|_\infty\leq \kappa_i$, pour tout $i\in\N$.
\end{itemize}
Enfin, pour une telle fibration $\pi'$, les parties 3 et 4 de 
\cite{cfg92} donnent la construction d'une métrique $g_2$ sur $M$ qui
est $T^k$-invariante et telle que $|\mathrm{D}^i(g_2-g_1)|\leq c(n,A,(N,h),i)$,
pour tout $i\in\N$. Cette dernière égalité assure l'existence d'une
constante $\tau_2(n,A,(N,h))$ telle que 
$\frac1{\tau_2}g_1\leq g_2\leq\tau_2g_1$, et permet aussi de contrôler la
courbure pour la métrique $g_2$.

On peut alors appliquer la proposition \ref{geom:pr}, qui nous
donne une métrique $g_3$ qui vérifie les points (1) à (3)
de la définition~\ref{geom:df}. On utilise le fait que comme $M$ est
proche de $N$ pour la distance de Gromov-Hausdorff, son diamètre
est contrôlé par $\varepsilon$ et le diamètre de $N$.

Pour obtenir la métrique $\tilde g$ du théorème \ref{geom:th}, il reste
à modifier la distribution horizontale de manière à ce que 
$\tilde g$ vérifie le point (4) de la définition~\ref{geom:df}.
Remarquons tout d'abord que comme l'application $e:\mathcal G^*\rightarrow
\mathcal H^2(N)$ est linéaire, il suffit de montrer l'égalité
$\de\omega=\pi'^*(e(\omega))$ pour une base de $\mathcal G^*$.
Soit $(\omega_i)$ une base de $\mathcal G^*$ orthonormée pour la
métrique $g_3$. Pour chaque $i$, $\de\omega_i$ s'écrit
\begin{equation}
\de\omega_i=\pi'^*(\alpha_i+\de\beta_i),
\end{equation}
où $\alpha_i$ est une forme harmonique et $\beta_i$ une forme cofermée.
On définit une nouvelle forme verticale $\omega_i'=\omega_i-\pi'^*(\beta_i)$.
Cette forme vérifie 
\begin{equation}\label{geom:eq3}
\de\omega_i'=\de\omega_i-\pi'^*(\de\beta_i)=\pi'^*(\alpha_i)\in\mathcal H^2(N).
\end{equation}
L'intersection des noyaux des formes $\omega_i'$ définit une nouvelle
distribution horizontale. On définit $\tilde g$ comme étant la métrique
sur $M$ telle que $\pi':(M,\tilde g)\rightarrow(N,h)$ soit une
submersion riemannienne et $\tilde g=g_3$ sur l'espace vertical.
Cette métrique vérifie le point (4) de la définition \ref{geom:df}
du fait de (\ref{geom:eq3}).

On doit encore vérifier $\tilde g$ est proche de $g_3$. Remarquons
que 
$$\tilde g-g_3=\sum_i(\omega'^2-\omega^2)=\sum_i(2\pi'^*(\beta_i)
\otimes\omega_i+\pi'^*(\beta_i)^2).$$
Or, B.~Colbois et G.~Courtois ont montré dans \cite{cc00} (lemme
A.32) qu'il existe
une constante $\kappa(n,a,(N,h))>0$ telle que les formes $\beta_i$ telles
qu'on les a définies vérifient $\|\beta_i\|_\infty\leq\kappa$, ce qui permet de
conclure. Remarque: le lemme A.32 de \cite{cc00} utilise le fait
que pour la métrique $g_3$, la norme de la seconde forme fondamentale 
est contrôlée et que la submersion $\pi'$ est $\kappa_i$-régulière. Il n'est
donc pas évident qu'on puisse obtenir le théorème \ref{geom:th} en supposant
que la métrique initiale $g$ est invariante et en se passant des 
résultats de \cite{cfg92}.

Enfin, il reste à montrer que la courbure de $(M,\tilde g)$ reste
bornée dans la direction horizontale. 
Soit $x\in N$, $\tilde X$ et $\tilde Y$ deux vecteurs orthonormés
tangents à $N$ en $x$, $\bar\omega$ une 
$1$-forme invariante de $T^k$, $\omega$ la $1$-forme induite sur $M$
pour la distribution horizontale associée à $g$ et $\omega'$
la $1$-forme induite pour la distribution associée à $\tilde g$.
Ces deux formes vérifient $\de\omega'=\pi'^*(\alpha)$ et
$\de\omega=\pi'^*(\alpha+\de\beta)$, où $\alpha$ est une $2$-forme
harmonique de $N$ et $\beta$ une $1$-forme de $N$.

D'après la formule de O'Neill, il suffit pour contrôler la courbure
sectionnelle $K_{(M,\tilde g)}(\pi'^*(X),\pi'^*(Y))$ de majorer
la norme de $[\pi'^*(X),\pi'^*(Y)]^V$. Or, on peut écrire d'une part,
\begin{equation}
\omega'([\pi'^*(X),\pi'^*(Y)]^V)=\de\omega'(\pi'^*(X),\pi'^*(Y))=\alpha(X,Y).
\end{equation}
D'autre part, on a
\begin{equation}\label{geom:li}
\|\alpha\|_\infty\leq \tau'(N,h)\|\alpha\|_2,
\end{equation}
d'après \cite{li80} (théorème 7), car $\alpha$ est harmonique, et
\begin{equation}
\|\alpha\|_2\leq\|\alpha+\de\beta\|_2=\|\de\omega\|_2\leq
\|\de\omega\|_\infty,
\end{equation}
en utilisant le fait qu'une forme harmonique est le plus petit
élément de sa classe de cohomologie pour la norme $L^2$.
Enfin, le lemme \ref{geom:oneill} permet de contrôler
la norme de $\de\omega$ en fonction de $a$ et $\|\omega\|_\infty$,
et la norme de $\omega$ est contrôlé en fonction de $\|\omega'\|_\infty$.
Comme la majoration
de $\omega'([\pi'^*(X),\pi'^*(Y)]^V)$ obtenue est indépendante du choix
de $\bar\omega$, on a bien une majoration de $|[\pi'^*(X),\pi'^*(Y)]^V|$
en fonction de $n$, $a$, $\varepsilon$ et $(N,h)$.

\section{Formes invariantes et petites valeurs propres}\label{inv}

Nous allons ici démontrer que pour étudier les petites valeurs propres
d'un fibré en tores muni d'une métrique invariante, on peut
se restreindre aux formes invariantes :
\begin{theo}\label{inv:th2}
Soit $k\in\N^*$, $T^k\hookrightarrow M\stackrel{\pi}{\rightarrow}N$ un
fibré en tore $T^k$, $\bar g$ une métrique invariante sur $T^k$ et $f$
une fonction sur $N$ strictement positive.
Supposons que $M$ est muni d'une métrique
$T^k$-invariante $g$ telle que pour tout $x\in N$, la restriction
$\bar g_x$ de $g$ à la fibre $\pi^{-1}(x)$ vérifie
$\bar g_x\leq f(x)\cdot\bar g$.

Soit $\lambda$ une valeur propre du laplacien agissant sur les formes
différentielles de $M$. Si $\lambda<\displaystyle(\sup_{x\in N}f(x))^{-1}
\cdot\lambda_{0,1}(T^k,\bar g)$, alors les formes propres
associées sont $T^k$-invariantes.
\end{theo}
\begin{rema}
On peut montrer que cette estimation est optimale : si on considère
un fibré trivial muni d'une métrique produit, on voit que les formes
propres de $T^k$ de valeur propre $\lambda_{0,1}(T^k,\bar g)$ induisent
sur le fibré des formes propres de même valeur propre qui ne sont pas
invariantes.
\end{rema}
\begin{rema}\label{inv:rem}
La démonstration du théorème met en évidence le fait que si la
multiplicité d'une valeur propre est impaire, alors le sous-espace
propre associé contient des formes invariantes.
\end{rema}
\begin{rema}
Le résultat peut s'étendre à toute  action isométrique d'un groupe
de Lie compact.
\end{rema}
\begin{proof}
Soit $E_\lambda$ un espace propre du laplacien de Hodge-de~Rham sur $M$.
Comme $M$ est munie d'une action isométrique de $T^k$, on peut décomposer
$E_\lambda$ en somme de représentations irréductibles $T^k$. En notant
$T^k=\R^k/\Gamma$, ces représentations sont classifiées par le réseau
dual $\Gamma^*$ de $\Gamma$. On a donc $E_\lambda=\oplus_{\gamma\in\Gamma^*}
E_{\lambda,\gamma}$, avec 
\begin{equation}
E_{\lambda,\gamma}=\{\omega\in E_\lambda,\ \forall\xi\in\R^k,\ 
\forall t\in\R,\ \varphi_{\xi,t}^*\omega=
e^{i2\pi t\langle\gamma,\xi\rangle}\omega\},
\end{equation}
où $\varphi_{\xi,t}$ est le flot du champ de vecteur $X_\xi$ induit par
l'action infinitésimale de $\xi$. Cet espace s'écrit aussi
\begin{equation}\label{inv:lie}
E_{\lambda,\gamma}=\{\omega\in E_\lambda,\ \forall\xi\in\R^k,\
\mathcal L_{X_\xi}\omega=i2\pi \langle\gamma,\xi\rangle \omega\}.
\end{equation}
La formule de Cartan nous donne, pour $\omega\in E_{\lambda,\gamma}$,
\begin{equation}
\langle\mathcal L_{X_\xi}\omega,\omega\rangle =
\langle\de\upiota_{X_\xi}\omega+\upiota_{X_\xi}\de\omega,\omega\rangle
=\langle\upiota_{X_\xi}\omega,\updelta\omega\rangle+
\langle\de\omega,X_\xi^b\wedge\omega\rangle.
\end{equation}

En utilisant le fait que $|X_\xi^b\wedge\omega|^2+|\upiota_{X_\xi}\omega|^2
=|X_\xi|^2|\omega|^2$, l'inégalité de Cauchy-Schwarz donne
\begin{eqnarray}\label{inv:maj}
|\langle\upiota_{X_\xi}\omega,\updelta\omega\rangle+
\langle\de\omega,X_\xi^b\wedge\omega\rangle| & \leq &
\|X_\xi\|_\infty\|\omega\|_{L^2}(\|\updelta\omega\|_{L^2}^2+
\|\de\omega\|_{L^2}^2)^{1/2}\nonumber\\
& \leq &\|X_\xi\|_\infty\lambda^{1/2}\|\omega\|_{L^2}^2.
\end{eqnarray}
 En combinant (\ref{inv:lie}) et (\ref{inv:maj}), on obtient que
$2\pi \langle\gamma,\xi\rangle \leq\|X_\xi\|_\infty\lambda^{1/2}$.

 On peut fixer le réseau $\Gamma$ de sorte que la métrique quotient
sur $T^k$ soit la métrique $\bar g$. Le spectre de $(T^k,\bar g)$ (pour
les fonctions) est alors $\{4\pi^2|\gamma|^2,\ \gamma\in\Gamma^*\}$ (cf.
\cite{ghl87}). Pour un $\xi$ donné de norme~1, la norme de $X_\xi$ est~1
pour la métrique $\bar g$ donc, $\|X_\xi\|_\infty^2\leq \sup_{x\in N}f(x)$.
Si $\gamma\neq0$ et qu'on choisit $\xi$ de norme~1 tel que 
$\langle\gamma,\xi\rangle=|\gamma|$, on a alors
\begin{equation}
\lambda_{0,1}(T^k,\bar g)\leq4\pi^2|\gamma|^2\leq\sup_{x\in N}f(x)\lambda.
\end{equation}
Si $\lambda<\sup_{x\in N}f(x)^{-1}\cdot\lambda_{0,1}(T^k,\bar g)$, alors
$\gamma=0$ et $E_\lambda$ est constitué de formes invariantes.
\end{proof}

\begin{proof}[de la remarque \ref{inv:rem}]
 Il suffit de remarquer que dans la démonstration du théorème, 
si la dimension du sous-espace propre $E_\lambda$ est impaire 
la décomposition de cet espace représentations
irréductibles contient nécessairement la représentation triviale, et
donc $E_\lambda$ contient des formes invariantes. 
\end{proof}

\section{Petites valeurs propres et rayon d'injectivité}\label{inj}
 Nous allons maintenant démontrer le théorème \ref{intro:th5}. Le cas
$k=1$ découlant de \cite{cc00}, nous supposerons $k\geq2$.
Soit $N$ une variété telle que $b_2(N)\geq k$. On
se donne un fibré principal en tores $M$ sur $N$ tel que l'application
$e:\mathcal G^*\rightarrow H^2(N)$ définie dans la section~\ref{topo} soit
 injective (c'est possible car $b_2(N)\geq k$).

On se donne sur $M$ et $N$ des métriques $g$ et $h$ telles que
le couple $(g,h)$ soit adapté (au sens défini en \ref{geom:adap}) à
la fibration $M\rightarrow N$ et que la restriction $\bar g$
de la métrique $g$ à la fibre $T^k$ soit le quotient de 
la métrique canonique de $\R^k$ à $T^k=\R^k/\Z^k$ (on a vu dans
la section~\ref{geom} que le fait que les fibres soient totalement
géodésiques est équivalent au fait qu'elles soient isométriques entre elles).
On va construire la famille de métriques $(g_\varepsilon)$ sur $M$ en faisant
varier la métrique $g_\varepsilon$ le long de la fibre, la connexion et 
la composante horizontale de la métrique restant identiques à celles de $g$. 
 
Soit $y\in\R^{k-1}$ un $(k-1)$-uplet difficilement approchable,
c'est-à-dire tel qu'il existe une constante $c(y)$ telle que pour tout 
$p\in\Z^{k-1}$ et tout $q\in\Z$, on a
\begin{equation}\label{inj:da}
\|p-qy\|^{k-1}|q|\geq c(y),
\end{equation}
la norme considérée étant la norme euclidienne canonique. On sait
qu'il existe une infinité non dénombrable de tels $(k-1)$-uplet
(voir \cite{sc80}, p.~22 et p.~41--43). On se donne une base
orthonormée $(X_1,\ldots, X_k)$ de $\R^k$ telle que $X_1$ soit
colinéaire au vecteur $(1,y)$. Pour tout $\varepsilon$, on se donne 
la base $\mathcal B_\varepsilon=(\frac1\varepsilon X_1,X_2,\ldots, X_k)$
et on définit la métrique $\bar g_\varepsilon$ sur $\R^k$ en posant
que $\mathcal B_\varepsilon$ est orthonormée pour $\bar g_\varepsilon$,
c'est-à-dire qu'on obtient $\bar g_\varepsilon$ en contractant les
longueurs d'un rapport $\varepsilon$ dans la direction de $X_1$.
Le passage de $\bar g_\varepsilon$ au quotient $\R^k/\Z^k$ détermine
une métrique sur $T^k$, qui elle-même induit la métrique $(g_\varepsilon)$
sur $M$. On a en particulier $g_1=g$. Les résultats classiques sur les 
effondrements (\cite{ca84}, \cite{cg86} théorème~2.1) assurent que 
la courbure de $(M,g_\varepsilon)$ reste bornée quand $\varepsilon$ tend 
vers zéro.

On va montrer l'inégalité \ref{intro:eq3} du théorème \ref{intro:th5}
en estimant successivement le volume, le rayon d'injectivité et la
première valeur propre de $M$ en fonction de $\varepsilon$.

\begin{enonce}{Fait}
Pour tout $\varepsilon\in]0,1]$, on a $\Vol(M,g_\varepsilon)=
\Vol(N,h)\cdot\varepsilon$.
\end{enonce}
\begin{proof}
Il suffit de remarquer que $\Vol(M,g_\varepsilon)=\Vol(N,h)
\cdot\Vol(T^k,\bar g_\varepsilon)$, et que $\Vol(T^k,\bar g_\varepsilon)=
\Det_{\mathcal B_\varepsilon}\mathcal B_1=\varepsilon$.
\end{proof}

Pour obtenir le théorème \ref{intro:th5}, on aura à reparamétrer
$g_\varepsilon$ de sorte que $\Vol(M,g_\varepsilon)=\varepsilon$.

\begin{enonce}{Fait}\label{inj:inj}
Il existe des constantes $\varepsilon_0(y,k)>0$ et $C(y,k)>0$ telle que 
$\injrad(T^k,\bar g_\varepsilon)\geq C\cdot\varepsilon^\frac1k$ pour tout
$\varepsilon<\varepsilon_0$.
\end{enonce}
\begin{proof} Calculer le rayon d'injectivité de 
$(T^k,\bar g_\varepsilon)$ revient à calculer le minimum des normes
dans $\R^k$ des points non nuls du réseau $\Z^k$ pour la métrique 
$\bar g_\varepsilon$. Plus précisément, on peut écrire
\begin{equation}\label{inj:min}
2\injrad(T^k,\bar g_\varepsilon)=\min_{(p,q)\in \Z^{k-1}\times\Z
\backslash(0,0)}\|(p,q)\|_{\bar g_\varepsilon}.
\end{equation}

Soit $(p,q)\in\Z^{k-1}\times\Z\backslash(0,0)$. On note $z$ la projection
orthogonale de $(p,q)$ sur la droite vectorielle $D$ engendrée par $X_1$ 
(remarque : elle ne dépend pas de $\varepsilon$), $\theta$ l'angle
entre $D$ et l'hyperplan $\R^{k-1}\times\{0\}$ pour la métrique canonique
(on a $\tan\theta=\|y\|^{-1}$) et $\theta'$ l'angle entre $D$ et
le vecteur $(p,q)-(q\cdot y,q)$ (voir figure \ref{inj:fig}). De plus, on
notera $\|\cdot\|_{\bar g_\varepsilon}$ la norme pour la métrique
$g_\varepsilon$ et $\|\cdot\|$ la norme euclidienne canonique. 
On peut écrire :
\begin{equation}
\|(p,q)\|^2_{\bar g_\varepsilon}=\|(p,q)-z\|^2_{\bar g_\varepsilon}
+\|z\|^2_{\bar g_\varepsilon}.
\end{equation}

\begin{figure}[h]
\begin{center}
\begin{picture}(0,0)%
\includegraphics{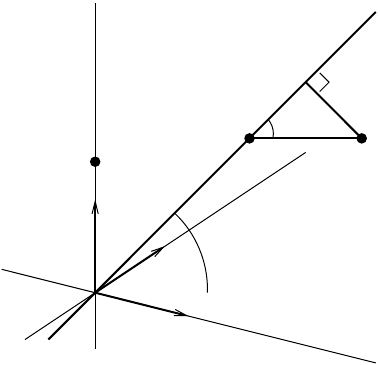}%
\end{picture}%
\setlength{\unitlength}{4144sp}%
\begingroup\makeatletter\ifx\SetFigFont\undefined%
\gdef\SetFigFont#1#2#3#4#5{%
  \reset@font\fontsize{#1}{#2pt}%
  \fontfamily{#3}\fontseries{#4}\fontshape{#5}%
  \selectfont}%
\fi\endgroup%
\begin{picture}(2885,2769)(1789,-6328)
\put(2086,-4819){\makebox(0,0)[lb]{\smash{{\SetFigFont{10}{12.0}{\rmdefault}{\mddefault}{\updefault}{\color[rgb]{0,0,0}$(0,q)$}%
}}}}
\put(3441,-5567){\makebox(0,0)[lb]{\smash{{\SetFigFont{10}{12.0}{\rmdefault}{\mddefault}{\updefault}{\color[rgb]{0,0,0}$\theta$}%
}}}}
\put(4332,-3749){\makebox(0,0)[lb]{\smash{{\SetFigFont{10}{12.0}{\rmdefault}{\mddefault}{\updefault}{\color[rgb]{0,0,0}$D$}%
}}}}
\put(3939,-4141){\makebox(0,0)[lb]{\smash{{\SetFigFont{10}{12.0}{\rmdefault}{\mddefault}{\updefault}{\color[rgb]{0,0,0}$z$}%
}}}}
\put(3155,-4605){\makebox(0,0)[lb]{\smash{{\SetFigFont{10}{12.0}{\rmdefault}{\mddefault}{\updefault}{\color[rgb]{0,0,0}$(qy,q)$}%
}}}}
\put(4617,-4676){\makebox(0,0)[lb]{\smash{{\SetFigFont{10}{12.0}{\rmdefault}{\mddefault}{\updefault}{\color[rgb]{0,0,0}$(p,q)$}%
}}}}
\put(4046,-5924){\makebox(0,0)[lb]{\smash{{\SetFigFont{10}{12.0}{\rmdefault}{\mddefault}{\updefault}{\color[rgb]{0,0,0}$\R^{k-1}$}%
}}}}
\put(3904,-4569){\makebox(0,0)[lb]{\smash{{\SetFigFont{10}{12.0}{\rmdefault}{\mddefault}{\updefault}{\color[rgb]{0,0,0}$\theta'$}%
}}}}
\put(2264,-3714){\makebox(0,0)[lb]{\smash{{\SetFigFont{10}{12.0}{\rmdefault}{\mddefault}{\updefault}{\color[rgb]{0,0,0}$\R$}%
}}}}
\end{picture}%
\end{center}
\caption{\label{inj:fig}}
\end{figure}

 Le vecteur $(p,q)-z$ est orthogonal à $D$, donc la norme 
$\|(p,q)-z\|_{\bar g_\varepsilon}$ est indépendante de $\varepsilon$. 
On a donc, en utilisant la relation (\ref{inj:da}):
\begin{equation}
\|(p,q)-z\|_{\bar g_\varepsilon}=\|(p,q)-z\|=\sin\theta'\|p-qy\|\geq
\frac{\sin\theta'\cdot c(y)^{\frac1{k-1}}}{|q|^{\frac1{k-1}}}.
\end{equation}
Comme $\theta\leq\theta'\leq\frac\pi2$ on en déduit
\begin{equation}
\|(p,q)-z\|^2_{\bar g_\varepsilon}\geq\frac{a}{|q|^{\frac2{k-1}}}
\end{equation}
où $a$ est une constante ne dépendant que de $k$ et $y$.

Le vecteur $z$ est parallèle à $D$, donc $\|z\|_{\bar g_\varepsilon}=
\varepsilon\|z\|$. De plus, on a 
\begin{equation}
\|z\|\geq\|(q,qy)\|-\|z-(q,qy)\|\geq\frac{\|qy\|}{\cos\theta}
-\|z-(q,qy)\|.
\end{equation}
On peut supposer que $\|z-(q,qy)\|\leq\frac{\|qy\|}{2\cos\theta}$. 
En effet, dans le cas contraire, on a
$\|(p,q)-z\|_{\bar g_\varepsilon}\geq\sin\theta\|p-qy\|\geq\frac
{\sin\theta}{2\cos\theta}\|qy\|\geq\frac q2$, ce qui signifie si $q\neq0$
que $\|(p,q)\|^2_{\bar g_\varepsilon}$ est minoré par $\frac12$ 
et donc ne réalise pas le minimum en (\ref{inj:min}) quand $\varepsilon$
est suffisamment petit (si $q=0$, c'est $\sin\theta\|p-qy\|$ qui est
uniformément minoré). On obtient finalement la minoration
\begin{equation}
\|z\|^2_{\bar g_\varepsilon}\geq\varepsilon^2\frac{\|qy\|^2}{4\cos^2\theta}
=\varepsilon^2bq^2,
\end{equation}
où $b$ est une constante qui ne dépend que de $y$, et donc
\begin{equation}
\|(p,q)\|^2_{\bar g_\varepsilon}\geq\frac{a}{|q|^{\frac2{k-1}}}
+\varepsilon^2bq^2.
\end{equation}

 La fonction $f:t\rightarrow at^{-\frac2{k-1}}+\varepsilon^2bt^2$
admet sur $\R^*_+$ un unique minimum en $t=\left(\frac a{\varepsilon^2b(k-1)}
\right)^{\frac{k-1}{2k}}$ dont la valeur est $\varepsilon^2a^{\frac{k-1}k}
b^\frac1k((k-1)^\frac1k+(k-1)^\frac{1-k}k)$. On peut donc en déduire que
\begin{equation}
\|(p,q)\|_{\bar g_\varepsilon}\geq C\varepsilon^\frac1k
\end{equation}
où $C$ est une constante ne dépendant que de $y$ et $k$.
\end{proof}
\begin{rema} Il découle du second théorème du corps convexe de
Minkowski (voir ch.~4 de \cite{sc80}, en particulier le lemme~1D) qu'il
existe une constante $c>0$ telle que pour toute métrique plate sur 
$T^k$, on a $\injrad(T^k)^{k-1}\cdot\diam(T^k)\leq c\cdot\Vol(T^k)$.
On peut donc déduire du fait~\ref{inj:inj} que $\diam(T^k,\bar g_\varepsilon)
=O(\varepsilon^{1/k})$.
\end{rema}
\begin{rema}
 Comme la fibre est totalement géodésique, son rayon d'injectivité
est égal au rayon d'injectivité de $M$ si $\varepsilon$ est suffisamment
petit.
\end{rema}

\begin{enonce}{fait}\label{inj:fait}
Il existe une constante $C'(N,k,e)>0$ telle que pour tout $\varepsilon$,
on a $\lambda_{p,1}(M,g_\varepsilon)\leq C'\cdot\varepsilon^2$ pour 
$p=1$ et $2$.
\end{enonce}
\begin{proof}
Commençons par le cas $p=1$. Pour estimer la première valeur propre 
non nulle, on va d'abord
calculer quelles sont les $1$-formes harmoniques. On sait déjà d'après
les résultats de la section précédente que les formes harmoniques
sont $T^k$-invariantes.
Soit $\varphi$ une $1$-forme différentielle $T^k$-invariante 
de $M$. On peut écrire 
\begin{equation}\label{inj:phi}
\varphi=\pi^*(\alpha)+\sum_{i=1}^k \pi^*(a_i)\cdot\omega_i,
\end{equation}
où $\alpha$ est une $1$-forme de $N$, $a_i$ des fonctions de $N$ et
$\omega_i$ les $1$-formes verticales induites par une base
orthonormée de $\mathcal G^*$. On a alors~:
\begin{equation}\label{inj:dphi}
\de\varphi=\pi^*(\de\alpha+\sum_{i=1}^k a_i\cdot e_i)+
\sum_{i=1}^k\de \pi^*(a_i)\wedge\omega_i,
\end{equation}
où $e_i$ désigne l'image de $\omega_i$ par l'application
$e:\mathcal G^*\rightarrow\mathcal H^2(N)$.
De plus, pour tout $i$ on a 
\begin{eqnarray}
\|\updelta(\pi^*(a_i)\omega_i)\|^2&=&
(\updelta(\pi^*(a_i)\omega_i),\updelta(\pi^*(a_i)\omega_i))\nonumber\\
&=&(\pi^*(a_i)\omega_i,\de\updelta(\pi^*(a_i)\omega_i)),
\end{eqnarray}
où $(\cdot,\cdot)$ désigne le produit scalaire $L^2$.
Comme $\pi^*(a_i)\omega_i$ est une forme $T^k$-invariante, 
$\updelta(\pi^*(a_i)\omega_i)$ est une fonction invariante, c'est-à-dire 
que c'est le relevé d'une fonction sur $N$. Par conséquent, 
$\de\updelta(\pi^*(a_i)\omega_i)$ est le relevé d'une $1$-forme
sur $N$, et est donc orthogonale à $\pi^*(a_i)\omega_i$. Finalement, on a:
\begin{equation}\label{inj:deltaphi}
\updelta\varphi=\pi^*(\updelta\alpha).
\end{equation}

Si $\varphi$ est harmonique, on a $\de\varphi=0$ et $\updelta\varphi=0$, donc 
\begin{equation}
\updelta\alpha=0,\ \de a_i=0\text{ pour tout $i$, et } 
\de\alpha+\sum_{i=1}^k a_i\cdot e_i=0.
\end{equation}
Comme les fonctions $a_i$ sont constantes, $\sum_{i=1}^k a_i\cdot e_i$ est
une $2$-forme harmonique de $N$, donc orthogonale à la forme exacte 
$\de\alpha$. On a donc $\Delta\alpha=0$ et $\sum_{i=1}^k a_i\cdot e_i=0$.
Comme $e$ est injective, les $e_i$ forment une famille libre, et donc
$a_i=0$ pour tout $i$. On obtient finalement que les $1$-formes harmoniques
de $M$ sont les relevés des $1$-formes harmoniques de $N$.

Pour majorer la première valeur propre non nulle du laplacien, il
suffit de calculer le quotient de Rayleigh pour la métrique $g_\varepsilon$
d'une $1$-forme orthogonale aux formes harmoniques. En notant $(\omega_i)$
les $1$-formes verticales induites par la base duale de 
$\mathcal B_\varepsilon$, on choisit comme forme 
test $\varepsilon^{-1}\omega_1$.  On vient de voir que la codifférentielle 
d'une telle forme est nulle, donc on peut écrire
\begin{equation}
R(\varepsilon^{-1}\omega_1)=\frac{\|\de(\varepsilon^{-1}\omega_1)\|^2
_{g_\varepsilon}}{\|\varepsilon^{-1}\omega_1\|^2_{g_\varepsilon}}.
\end{equation}
La forme $\varepsilon^{-1}\omega_1$ est indépendante de $\varepsilon$,
donc $\de(\varepsilon^{-1}\omega_1)$ aussi, et comme de plus elle est 
horizontale donc sa norme ponctuelle ne dépend pas de $\varepsilon$; on a 
donc $\|\de(\varepsilon^{-1}\omega_1)\|^2=c\Vol(M,g_\varepsilon)$. Par ailleurs,
 $\|\varepsilon^{-1}\omega_1\|^2_{g_\varepsilon}=\varepsilon^{-1}
\|\omega_1\|^2_{g_\varepsilon}$, donc 
$\|\varepsilon^{-1}\omega_1\|^2_{g_\varepsilon}$ est de la forme 
$c'\varepsilon^{-2}\Vol(M,g_\varepsilon)$. Comme 
$\lambda_{1,1}(M,g_\varepsilon)\leq R(\varepsilon^{-1}\omega_1)$, on en 
déduit que $\lambda_{1,1}(M,g_\varepsilon)\leq C'\cdot\varepsilon^2$, où
$C'$ est une constante dépendant des choix de $N$, $k$ et $e$.

 Le cas $p=2$ se déduit du premier par la théorie de Hodge, en remarquant
que dans le cas $p=1$ on a trouvé une petite valeur propre du laplacien
retreint aux formes cofermées.
\end{proof}

\begin{enonce}{fait}
Si $b_1(N)>b_2(M)$, alors il existe une constante $C''(N,k,e)>0$ telle que 
pour tout $\varepsilon$, on a $\lambda_{p,b_1(N)-b_2(M)}(M,g_\varepsilon)
\leq C''\cdot\varepsilon^2$ pour $p=2$ et $3$.
\end{enonce}

\begin{proof}
On va se restreindre une nouvelle fois aux formes cofermées pour
montrer le résultat pour $p=2$ et ensuite en déduire le cas $p=3$.

Contrairement au cas des $1$-formes, on a pas en général de moyen simple
de déterminer quelles sont les $2$-formes harmoniques (on rencontre
la même difficulté dans l'étude du spectre des fibrés en cercles, voir
\cite{cc00}). On va majorer
le quotient de Rayleigh sur un espace test de dimension 
$b_1(N)$, et la condition $b_1(N)>b_2(M)$ assurera qu'on a bien majoré
une ou plusieurs valeur propres non nulles.

 En notant toujours $(\omega_i)$ les $1$-formes verticales induites 
par la base duale de $\mathcal B_\varepsilon$ on considère une forme 
différentielle $\varphi=\omega_1\wedge\pi^*(\alpha)$, où $\alpha$ est une 
$1$-forme harmonique de $N$. Son quotient de Rayleigh est
\begin{equation}
R(\varphi)=\frac{\|\de(\omega_1\wedge\pi^*(\alpha))\|^2
_{g_\varepsilon}+\|\updelta(\omega_1\wedge\pi^*(\alpha))\|^2_{g_\varepsilon}}
{\|\omega_1\wedge\pi^*(\alpha)\|^2_{g_\varepsilon}}
\end{equation}
La forme $\varphi$ est cofermée. En effet, si on calcule la différentielle
de sa forme duale, on obtient :
\begin{eqnarray}
\de *(\omega_1\wedge\pi^*(\alpha))&=&\de((\bigwedge_{i\neq1}\omega_i
)\wedge\pi^*(*\alpha))\nonumber\\
&=&\sum_{i\neq1}(-1)^i\bigwedge_{j\neq1,i}\omega_j\wedge\pi^*(e(\omega_i))
\wedge\pi^*(*\alpha)
\end{eqnarray}
Comme $\pi^*(e(\omega_i))\wedge\pi^*(*\alpha)=\pi^*(e(\omega_i)\wedge*\alpha)$
et que $e(\omega_i)\wedge*\alpha$ est de degré supérieur à la dimension de $N$,
on en déduit que tous les termes de la somme sont nuls.
On est ramené à :
\begin{equation}
R(\varphi)=\frac{\|\de(\omega_1\wedge\pi^*(\alpha))\|^2
_{g_\varepsilon}}{\|\omega_1\wedge\pi^*(\alpha)\|^2_{g_\varepsilon}}.
\end{equation}
On a d'une part $|\omega_1\wedge\pi^*(\alpha)|^2_{g_\varepsilon}
=|\omega_1|^2_{g_\varepsilon}|\pi^*(\alpha)|^2_{g_\varepsilon}=
|\alpha|^2$ et donc
\begin{equation}
\|\omega_1\wedge\pi^*(\alpha)\|^2_{g_\varepsilon}=\|\alpha\|^2,
\end{equation}
et d'autre part 
\begin{equation}
\|\de(\omega_1\wedge\pi^*(\alpha))\|^2_{g_\varepsilon}=
\|\de\omega_1\wedge\pi^*(\alpha)\|^2_{g_\varepsilon}
\leq\|\de\omega_1\|^2_{g_\varepsilon}\|\alpha\|^2.
\end{equation}
On voit finalement que $R(\varphi)\leq R(\omega_1)$ quel que soit $\alpha$.
La majoration de $R(\omega_1)$ obtenue précédemment permet de conclure.
\end{proof}

On va maintenant donner des exemples de fibrés pour lesquels la condition
$b_1(N)>b_2(M)$ est bien vérifiée.

\begin{exem}\label{inj:ex}
On considère pour tout entier $k$ la variété $M=(S^3)^k\times S^1$.
La fibration de Hopf $S^1\hookrightarrow S^3\rightarrow S^2$ induit une 
action libre de $T^k=(S^1)^k$ sur $(S^3)^k$, et permet donc de définir
une fibration principale $T^k\hookrightarrow M\rightarrow N$, où 
$N=(S^2)^k\times S^1$. La classe d'Euler $e$ de ce fibré est bien injective :
l'image de $e$ est engendrée par les relevés à $M$ de chacune des classes
d'Euler correspondant aux fibrations de Hopf qui sont bien linéairement
indépendants. Enfin, la formule de Künneth permet de calculer aisément
les nombres de Betti des produits $M$ et $N$, on obtient en particulier
$b_1(N)=1$ et $b_2(M)=0$.
\end{exem}

\begin{rema}
Dans l'exemple précédent, les métriques sur $M$ et $N$ ne sont 
\emph{a priori} pas des métriques produits. On ne peut donc pas déduire
directement les propriétés du spectre de la formule de Künneth.
\end{rema} 

Pour justifier la remarque~\ref{intro:rem} faite dans l'introduction,
on va construire des contre-exemples aux théorèmes~1.1 et~1.2 de~\cite{ct97}.
Ceux-ci permettent de montrer (cf.~\cite{ja05}) que si, pour un réel $r$ tel 
que $0<r<\injrad(M)$, on note~$q$ le nombre de boules géodésiques de rayon 
$4^{-n}r$ nécessaire pour recouvrir $M$, alors il existe une constante $C$ 
ne dépendant que de la dimension et de la borne sur la courbure telle que
\begin{equation}\label{inj:ct97}
\lambda_{p,1}(M,g)\geq Cr^{-2}q^{-7(p+1)}.
\end{equation}
C'est cette inégalité qu'on va contredire.
\begin{exem}\label{inj:ex2}
À courbure bornée et dimension $n$ fixée, le nombre $q$ de boules géodésiques 
de rayon $4^{-n}r$ nécessaire pour recouvrir $M$ est de l'ordre de 
$\Vol(M)\cdot r^{-n}$ quand $r$ est petit. Pour $p=1$, la minoration
de l'inégalité (\ref{inj:ct97}) est donc de l'ordre de 
$\Vol(M)^{-14}r^{14n-2}$.

On se donne un fibré $T^k\hookrightarrow M\to N$ et un effondrement 
$g_\varepsilon$ fournis par le théorème~\ref{intro:th5}. On a alors
$\Vol(M)=\varepsilon$ et on peut choisir $r\sim\varepsilon^{\frac1k}$,
ce qui donne $q\sim\varepsilon^{1-n/k}$ et une minoration de l'ordre de 
$\varepsilon^{(14n-14k-2)/k}=\varepsilon^{(14d-2)/k}$ en notant 
$d=n-k$ la dimension de la base $N$.

 Or, si on fixe la dimension $d\geq3$ de $N$ et qu'on choisit $k$ suffisamment 
grand (c'est possible puisque $b_2(N)$ peut être arbitrairement grand à 
dimension fixée), l'estimation $\varepsilon^{(14d-2)/k}$ ne peut pas minorer
$\lambda_{1,1}(M,g_\varepsilon)$ qui est de l'ordre de $\varepsilon^2$ quand
$\varepsilon$ tend vers~0.

\end{exem}

Il reste enfin à démontrer le corollaire \ref{intro:cor}.\\
\begin{proof}
 Cas $n=5$ : On considère une variété compacte $N$ de dimension 3 telle 
que $b_2(N)\geq2$ et on choisit pour $M$ un fibré principal en tore
$T^2$ sur $N$ dont la classe d'Euler est injective. Les résultats
précédents assure l'existence d'un effondrement tel que 
$\lambda_{p,1}(M)\leq C\cdot\injrad(M)^4$ pour $p$ égal à 1 et 2. Par dualité
de Hodge, on a la même inégalité pour $n-1$ et $n-2$, donc pour tout $p$
entre 1 et $n-1$.

Cas $n=7$ : On peut vérifier que le fibré construit dans l'exemple
\ref{inj:ex} pour $k=2$ convient. On sait déjà que 
$\lambda_{p,1}(M)\leq C\cdot\injrad(M)^4$ pour $p=1$, 2 et 3, et la dualité 
de Hodge donne l'inégalité pour $p=4$, 5 et 6.
\end{proof}

\section{Minoration du spectre des fibrés principaux en tores}
\label{adap}
\subsection{Minoration de la première valeur propre des $1$-formes}
 Les résultats des sections \ref{topo} à \ref{inv} nous permettent 
de démontrer le théorème~\ref{intro:th4}. On a vu qu'on pouvait
se ramener au cas d'un fibré muni d'une métrique adaptée. On va donc
montrer le résultat du théorème \ref{intro:th4} pour un fibré vérifiant
les conclusions du théorème \ref{geom:th}:
\begin{theo}\label{adap:th}
Soit $a>0$, $d>0$ deux réels, $n$, $k$ et $m$ trois entiers
tels que $n=k+m$, et $(N^m,h)$ une variété riemannienne. 
Il existe des constantes 
$c(n,a,d,(N,h))$
et $\varepsilon(n,a,d,(N,h))$ strictement
positives telles que si $\bar g$ est une métrique sur le tore 
$T^k$  telle que $\diam(T^k)<\varepsilon$ et si 
$T^k\hookrightarrow M^n\rightarrow N$ est un fibré principal
muni d'un couple de métriques $(g,h)$ adapté au fibré et tel que 
$g=\bar g$ en restriction à la fibre, 
$\diam(M,g)<d$ et $|K_M(X,Y)|\leq a$ pour toute
paire $(X,Y)$ de vecteurs horizontaux orthonormés, alors on a
$$\lambda_{1,1}(M,g)\geq c\cdot \Vol^2(T^k).$$
\end{theo}
On s'est ici donné comme hypothèse que la métrique sur $N$ est fixée.
En effet, une hypothèse sur la courbure ne nous sera pas suffisante.
Par ailleurs
la minoration du spectre en fonction du rayon d'injectivité 
découle du fait que le rayon d'injectivité de $M$ est égal à celui de 
la fibre si $\varepsilon$ est suffisamment petit (la fibre étant
totalement géodésique) et que, la restriction de la métrique 
à la fibre étant plate, on a $\Vol(T^k)\geq C\injrad(T^k)^k$ où
$C$ ne dépend que de $k$.

Dans un premier temps, nous allons démontrer le théorème dans
le cas où le fibré $M$ ne contient pas de sous-fibré trivial,
c'est-à-dire quand  l'application $e:\mathcal G^*\rightarrow\mathcal H^2
(N,h)$ est injective. Nous généraliserons ensuite le résultat 
à un fibré principal quelconque. D'autre part, on se restreindra 
aux formes $T^k$-invariantes, en vertu des résultats de la section~\ref{inv} 
(corollaire \ref{inv:th2}). En effet, le spectre des formes 
orthogonales aux formes invariantes sera minoré en fonction de 
la constante $\varepsilon$ du théorème, et on pourra toujours
choisir cette constante suffisamment petite de sorte que le
spectre des formes orthogonales aux formes invariantes soit
plus grand que le terme $c\cdot \Vol^2(T^k)$.

 Supposons donc $e$ injective. La démonstration se déroule en
deux étapes. D'abord, on se ramène à l'étude des valeurs
propres de l'opérateur $e^*e$, l'adjoint étant défini en munissant
$\mathcal H^2(N)$ de sa norme $L^2$  :

\begin{enonce}{fait}\label{adap:ft1}
Il existe $\varepsilon(n,a,\lambda_{0,1}(N,h),\lambda_{1,1}(N,h))>0$ et 
$c(n,a,\lambda_{0,1}(N,h))>0$ tels que
pour toute $1$-forme $\varphi$ sur $M$ $T^k$-invariante et
orthogonale à $\Ker\Delta^1(M,g)$,
si le quotient de  Rayleigh de $\varphi$ vérifie 
$R(\varphi)<\varepsilon$, alors il existe une forme $\omega$
induite par un élément de $\mathcal G^*$ telle que 
$\|e(\omega)\|^2\leq c\cdot\varepsilon\|\omega\|^2$.
\end{enonce}

\begin{proof}
Soit $\varphi$ une $1$-forme différentielle $T^k$-invariante 
de $M$. On peut comme dans la section précédente écrire $\varphi$ sous
la forme
\begin{equation}\label{adap:phi}
\varphi=\pi^*(\alpha)+\sum_{i=1}^k \pi^*(a_i)\cdot\omega_i,
\end{equation}
où $\alpha$ est une $1$-forme de $N$, $a_i$ des fonctions de $N$ et
$\omega_i$ les $1$-formes verticales induites par une base
orthonormée de $\mathcal G^*$. On a alors~:
\begin{equation}\label{adap:dphi}
\de\varphi=\pi^*(\de\alpha+\sum_{i=1}^k a_i\cdot e_i)+
\sum_{i=1}^k\de \pi^*(a_i)\wedge\omega_i,
\end{equation}
où $e_i$ désigne l'image de $\omega_i$ par l'application
$e:\mathcal G^*\rightarrow\mathcal H^2(N)$. On a vu de plus dans la
démonstration du fait \ref{inj:fait} que $\varphi$ vérifie
\begin{equation}\label{adap:deltaphi}
\delta\varphi=\pi^*(\delta\alpha),
\end{equation}
et que le fait que $e$ soit injectif implique que les $1$-formes harmoniques
de $M$ sont exactement les relevés des $1$-formes harmoniques de $N$.

Supposons que $\varphi$ est de norme $1$, c'est-à-dire que 
$\|\alpha\|^2+\sum_{i=1}^k\|a_i\|^2=1$, et qu'elle est orthogonale aux 
formes harmoniques de $M$. Le quotient de 
Rayleigh de $\varphi$ s'écrit alors
\begin{equation}\label{adap:eq1}
R(\varphi)=\|\updelta\alpha\|^2+\|\de\alpha+\sum_{i=1}^k a_i\cdot e_i\|^2
+\sum_{i=1}^k\|\de \pi^*(a_i)\|^2.
\end{equation}
 Supposons que $R(\varphi)<\varepsilon$ pour un $\varepsilon>0$ donné.
On a en particulier $\|\de\alpha+\sum_{i=1}^k a_i\cdot e_i\|^2<\varepsilon$.
Pour tout $i$, notons $\bar a_i$ la valeur moyenne de la fonction $a_i$.
On peut alors écrire
\begin{equation}
\|\de\alpha+\sum_{i=1}^k\bar a_i\cdot e_i
+\sum_{i=1}^k(a_i-\bar a_i)\cdot e_i\|^2<\varepsilon.
\end{equation}
L'inégalité triangulaire nous donne alors :
\begin{equation}
\|\de\alpha+\sum_{i=1}^k\bar a_i\cdot e_i\|<
\|\sum_{i=1}^k(a_i-\bar a_i)\cdot e_i\|+\sqrt\varepsilon.
\end{equation}
Comme les fonctions $(a_i-\bar a_i)$ sont de moyennes nulles, leur quotient
de Rayleigh est supérieur à $\lambda_{0,1}(N,h)$, et donc 
$\|a_i-\bar a_i\|^2\leq\frac{\|\de a_i\|^2}{\lambda_{0,1}(N,h)}\leq
\frac\varepsilon{\lambda_{0,1}(N,h)}$.
Par ailleurs, la forme exacte $\de\alpha$ est orthogonale à la
forme harmonique $\sum_{i=1}^k\bar a_i\cdot e_i$. On a donc finalement
\begin{eqnarray}\label{adap:eq2}
\|\sum_{i=1}^k\bar a_i\cdot e_i\|^2 & \leq & 
\|\de\alpha+\sum_{i=1}^k\bar a_i\cdot e_i\|^2\nonumber\\
&\leq&(\|\sum_{i=1}^k(a_i-\bar a_i)\cdot e_i\|+\sqrt\varepsilon)^2\nonumber\\
&\leq&(\sum_{i=1}^k(\|a_i-\bar a_i\|\cdot\|e_i\|_\infty)
+\sqrt\varepsilon)^2\nonumber\\
& \leq & \varepsilon\left(1+
\frac{\sum_{i=1}^k\|e_i\|_\infty}{\sqrt{\lambda_{0,1}(N,h)}}\right)^2.
\end{eqnarray}
Le lemme \ref{geom:oneill} permet d'obtenir une majoration de 
$\|\sum_{i=1}^k\bar a_i\cdot e_i\|^2$ en fonction de $n$, $\lambda_{0,1}(N,h)$
et d'une borne sur la courbure de $(M,g)$.

Pour tout $i$, la fonction $a_i-\bar a_i$ est orthogonale à $\bar a_i$,
donc 
\begin{equation}\label{adap:eq3}
\|\alpha\|^2+\sum_{i=1}^k(\|\bar a_i\|^2+\|a_i-\bar a_i\|^2)=1.
\end{equation}
D'une part, on a déjà vu que chaque terme $\|a_i-\bar a_i\|^2$ est 
majoré par $\frac\varepsilon{\lambda_{0,1}(N,h)}$. D'autre part, comme
$\varphi$ est orthogonale aux relevés des formes harmoniques de $(N,h)$, 
$\alpha$ est elle-même orthogonale aux formes harmoniques de $N$. On
peut donc écrire
\begin{equation}
\|\alpha\|^2\leq\frac1{\lambda_{1,1}(N,h)}(\|\de\alpha\|^2+\|\updelta\alpha\|^2).
\end{equation}
Le terme $\|\updelta\alpha\|^2$ est majoré par $R(\varphi)$, donc par 
$\varepsilon$, et $\|\de\alpha\|^2$ est majoré en fonction de $a$, $n$ et
$\lambda_{0,1}(N,h)$ comme dans l'inégalité \ref{adap:eq2}. Il découle
donc de l'équation (\ref{adap:eq3}) :
\begin{equation}\label{adap:eq4}
\sum_{i=1}^k\|\bar a_i\|^2\geq 
1-\tau(n,a,\lambda_{0,1}(N,h),\lambda_{1,1}(N,h))\varepsilon.
\end{equation}
 
Si on prend pour $\omega$ la $1$-forme $\sum_{i=1}^k\bar a_i\omega_i$,
on a 
\begin{equation}
\frac{\|e(\omega)\|^2}{\|\omega\|^2}=
\frac{\|\sum_{i=1}^k\bar a_i\cdot e_i\|^2}{\sum_{i=1}^k\|\bar a_i\|^2}
\end{equation}
 Selon (\ref{adap:eq4}), le dénominateur est supérieur à $\frac12$ si
$\varepsilon$ est suffisamment petit. La majoration du numérateur 
fournie par (\ref{adap:eq2}) donne alors le résultat souhaité.
\end{proof}

 On va maintenant minorer le spectre de $e^*e$ en fonction du
volume de la fibre $T^k$.

\begin{enonce}{fait}\label{adap:ft2}
Il existe une constante $c(n,a,(N,h))>0$ telle que
la première valeur propre de $e^*e$ soit minorée par 
$c\cdot\Vol(T^k)^2$.
\end{enonce}

\begin{proof}
Soient $\lambda_1,\cdots,\lambda_k$ les valeurs propres de $e^*e$
classées dans l'ordre croissant. Comme $e$ est injective, ces
valeurs propres sont non nulles et la premi\`ere vérifie
\begin{equation}
\lambda_1=\frac{\prod_i\lambda_i}{\prod_{i\neq1}\lambda_i}
=\frac{\Det(e^*e)}{\prod_{i\neq1}\lambda_i}.
\end{equation}
Par ailleurs, les valeurs propres de $e^*e$ v\'erifient
$\lambda_i\leq\|e^*e\|$, donc
\begin{equation}
\lambda_1\geq\frac{\Det(e^*e)}{\|e^*e\|^{k-1}}\geq\frac{\Det(e^*e)}
{\|e\|^{2k-2}}.
\end{equation}

L'image de $e$ est un sous-espace de $\Ker \Delta^2(N)$ de
dimension $k$, engendré par un sous-réseau du réseau des
formes harmoniques entières de $N$. Si on restreint $e$ et
$e^*$ à ce sous-espace, on peut écrire 
$\Det(e^*e)=(\Det e)^2$, où $\Det e$ est le déterminant d'une
matrice de $e$ écrite dans des bases orthonormées de 
$\mathcal G^*$ et $\Ima e$, ce qui donne 
\begin{equation}
\lambda_1\geq\frac{(\Det e)^2}{\|e\|^{2k-2}}.
\end{equation}

 Le lemme \ref{geom:oneill} donne une majoration de $\|e\|$ en fonction
de $n$ et de la borne $a$ sur la courbure de $M$, il ne reste 
donc qu'\`a minorer $\Det e$. Notons $\Det' e$ le déterminant de la
matrice de $e$ dans la base canonique de $\mathcal G^*={\R^k}^*$ 
et une base orthonormée de $\Ima e$. On a alors
$\Det e=(\Det' e)(\Vol T^k)$. Comme les images dans $\Ker\Delta^2(N)$
des éléments de la base canonique de $\mathcal G^*$ sont des formes 
entières, le déterminant $\Det' e$, qui est aussi le volume de 
$e([0,1]^k)$, est un multiple du volume d'un domaine fondamental 
du réseau des formes entières dans $\Ima e$. Comme par ailleurs
$\Det' e$ est non nul, il sera donc minoré par le volume
de ce domaine fondamental. Si on note $\rho$ le minimum des normes
des $2$-formes harmoniques entières non nulles, ce volume 
est minoré par le volume d'une boule de rayon $\frac\rho2$ dans 
$\Ima e$, et donc minoré par une constante ne dépendant que de 
$n$ et de la métrique $h$ de $N$. On peut donc bien écrire
$$\lambda_1\geq c(n,a,(N,h))\cdot\Vol(T^k)^2.$$
\end{proof}

 Nous allons maintenant supposer que $e$ n'est pas injective. Notons
$l$ la dimension de son noyau. Le premier nombre de Betti de $M$ est
alors $b_1(N)+l$. En effet, on a vu que si une $1$-forme $\varphi=
\pi^*(\alpha)+\sum_{i=1}^k \pi^*(a_i)\cdot\omega_i$ est harmonique, cela
signifie, d'après (\ref{adap:dphi}) et (\ref{adap:deltaphi}) :
\begin{equation}
\Delta\alpha=0,\de a_i=0\text{ pour tout $i$, et }\sum_{i=1}^k a_i\cdot e_i=0.
\end{equation}
 Comme les fonctions $a_i$ sont constantes, l'ensemble des $a_i$ tels
que
$\sum_{i=1}^k a_i\cdot e_i=0$ est exactement le noyau de $e$.
L'espace des formes harmoniques de $M$ est donc l'espace
engendré par les relevés des formes harmoniques de $N$ et les
formes verticales induites par les éléments de noyau de $e$.

 On peut reprendre la démonstration précédente en prenant pour
$(\omega_i)_i$ une base de $\mathcal G^*$ telle que $\omega_{k-l+1},
\cdots,\omega_k$ soit une base de $\Ker e$ (le fait que la
forme $\varphi$ est orthogonale aux formes harmoniques se traduit
par le fait que $a_{k-l+1},\cdots,a_k=0$) et en étudiant $e^*e$
restreint à l'orthogonal de $\Ker e$. On obtient de la même façon
le résultat du fait \ref{adap:ft1}, à savoir que la première valeur 
propre du laplacien sur $M$ est minorée à une constante multiplicative 
près par la première valeur propre de $(e^*e)_{|(\Ker e)^\perp}$.

Pour minorer le spectre de $(e^*e)_{|(\Ker e)^\perp}$, on doit
être un peu plus attentif dans la manipulation des bases de 
$\mathcal G^*$.

 Soit $\mathcal B=(\omega_1,\cdots,\omega_k)$ une base orthonormée de 
$\mathcal G^*$ et $\mathcal B'=(\omega_1',\cdots,\omega_k')$ une
base du réseau des entiers de $\mathcal G^*$, telles que 
$(\omega_1,\cdots,\omega_l)$ et $(\omega_1',\cdots,\omega_l')$ soient
des bases de $\Ker e$ (comme l'image du réseau des entiers de $\mathcal G^*$
est contenue dans un réseau, le noyau de $e$ est effectivement engendré
par des éléments entiers). La matrice de passage de $\mathcal B$ à 
$\mathcal B'$ est de la forme $$P=\left(\begin{array}{cc}P_1 & P_2\\
0 & P_3\end{array}\right),$$ où $P_1$ est un bloc carré de taille $l$.
Si se donne une base orthonormée de $\Ima e$, la matrice de $e$ 
s'écrit sous la forme $(0,A)$ dans la base $\mathcal B$ et 
$(0,A')$ dans la base $\mathcal B'$, où $A$ et $A'$ sont des
blocs carrés de taille $k-l$ et vérifient $A'=AP_3$.

Le spectre de $(e^*e)_{|(\Ker e)^\perp}$ est celui de $A^*A$. On peut 
écrire, comme dans la démonstration du fait \ref{adap:ft2}~:
\begin{equation}
\lambda_1\geq\frac{\Det A^*A}{\|A^*A\|^{k-l-1}}
\geq \frac{(\Det A)^2}{\|A\|^{2(k-l-1)}},
\end{equation}
où $\lambda_1$ est la première valeur propre non nulle de $e^*e$.
De plus, on a $\Det A'=\Det A\cdot\Det P_3$ et donc
\begin{equation}
\Det A=\frac{\Det A'}{\Det P_3}=\Det A'\frac{\Det P_1}{\Det P}.
\end{equation} 
Le déterminant de $A'$ est, comme précédemment, minoré par le covolume
du réseau des formes entières dans $\Ima e$, et $\Det P$ 
s'interprète géométriquement comme l'inverse du volume de $T^k$. 

Il reste à minorer $\Det P_1$. Comme $\Ker e$ est engendré par des
éléments entiers de $\mathcal G^*$, l'orthogonal de $\Ker e$ pour la 
dualité définit un sous-tore $T^{k-l}$ de $T^k$. De plus, le dual de 
de l'algèbre de Lie $\mathcal G(T^k/T^{k-l})$ du quotient $T^k/T^{k-l}$ 
est isomorphe à $\Ker e$. La matrice $P_1$ est donc la matrice de passage 
d'une base orthonormée de $\mathcal G^*(T^k/T^{k-l})$ dans une base du
réseau des entiers de $\mathcal G^*(T^k/T^{k-l})$, et par conséquent
$\Det P_1$ est l'inverse du volume de $T^k/T^{k-l}$ pour la métrique
quotient. Le diamètre de $T^k/T^{k-l}$ est majoré par $\varepsilon$,
comme celui de $T^k$, et par conséquent son volume aussi.

\bibliographystyle{smfalpha}
\bibliography{vol}

\providecommand{\bysame}{\leavevmode ---\ }
\providecommand{\og}{``}
\providecommand{\fg}{''}
\providecommand{\smfandname}{\&}
\providecommand{\smfedsname}{\'eds.}
\providecommand{\smfedname}{\'ed.}
\providecommand{\smfmastersthesisname}{M\'emoire}
\providecommand{\smfphdthesisname}{Th\`ese}
\begin{thebibliography}{GHL87}

\bibitem[ALK00]{alk00}
{\scshape J.~A. \'Alvarez~L\'opez {\normalfont \smfandname} Y.~A. Kordyukov} --
  {\og {A}diabatic limits and spectral sequences for riemannian foliations\fg},
  \emph{Geom. funct. anal.} \textbf{10} (2000), no.~5, p.~977--1027.

\bibitem[Bes08]{be08}
{\scshape A.~L. Besse} -- \emph{{E}instein {M}anifolds}, Classics in
  Mathematics, Springer Verlag, 2008.

\bibitem[BT82]{bt82}
{\scshape R.~Bott {\normalfont \smfandname} L.~W. Tu} -- \emph{{D}ifferential
  form in algebraic topology}, Springer Verlag, 1982.

\bibitem[Car84]{ca84}
{\scshape Y.~Carri\`ere} -- {\og {L}es propri\'et\'es topologiques des flots
  riemanniens retrouv\'ees \`a l'aide du th\'eor\`eme des vari\'et\'es presque
  plates\fg}, \emph{Math. Z.} \textbf{186} (1984), p.~393--400.

\bibitem[CC90]{cc90}
{\scshape B.~Colbois {\normalfont \smfandname} G.~Courtois} -- {\og A note on
  the first non zero eigenvalue of the {L}aplacian acting on $p$-forms\fg},
  \emph{Manuscripta Math.} \textbf{68} (1990), no.~2, p.~143--160.

\bibitem[CC00]{cc00}
\bysame , {\og Petites valeurs propres des $p$-formes diff\'erentielles et
  classe d'{E}uler des ${S^1}$-fibr\'es\fg}, \emph{Ann. scient. \'Ec. norm.
  sup. (4)} \textbf{33} (2000), no.~5, p.~611--645.

\bibitem[CFG92]{cfg92}
{\scshape J.~Cheeger, K.~Fukaya {\normalfont \smfandname} M.~Gromov} -- {\og
  {N}ilpotent structures and invariant metrics on collapsed manifolds\fg},
  \emph{J. Amer. Math. Soc.} \textbf{5} (1992), no.~2, p.~327--372.

\bibitem[CG86]{cg86}
{\scshape J.~Cheeger {\normalfont \smfandname} M.~Gromov} -- {\og {C}ollapsing
  riemannian manifolds while keeping their curvature bounded {I}\fg}, \emph{J.
  Differential Geom.} \textbf{23} (1986), no.~3, p.~309--346.

\bibitem[CT97]{ct97}
{\scshape S.~Chanillo {\normalfont \smfandname} F.~Tr\`eves} -- {\og On the
  lowest eigenvalue of the {H}odge {L}aplacian\fg}, \emph{J. Differ. Geom.}
  \textbf{45} (1997), no.~2, p.~273--287.

\bibitem[Dod82]{do82}
{\scshape J.~Dodziuk} -- {\og Eigenvalues of the {L}aplacian on forms\fg},
  \emph{Proc. of Am. Math. Soc.} \textbf{85} (1982), no.~3, p.~438--443.

\bibitem[For95]{fo95}
{\scshape R.~Forman} -- {\og Spectral {S}equences and {A}diabatic {L}imits\fg},
  \emph{Comm. Math. Phys.} \textbf{168} (1995), no.~1, p.~57--116.

\bibitem[GHL87]{ghl87}
{\scshape S.~Gallot, D.~Hulin {\normalfont \smfandname} J.~Lafontaine} --
  \emph{{R}iemannian {G}eometry}, Springer Verlag, 1987.

\bibitem[Gro80]{gr80}
{\scshape M.~Gromov} -- {\og Paul {L}evy's isoperimetric inegality\fg},
  \emph{Pr\'epublication IH\'ES} (1980), paru dans \cite{gr99}.

\bibitem[Gro99]{gr99}
\bysame , \emph{Metric {S}tructures for {R}iemannian and {N}on-{R}iemannian
  {S}paces}, Birkh\"auser, 1999.

\bibitem[Her60]{her60}
{\scshape R.~Hermann} -- {\og A sufficient condition that a mapping of
  {R}iemannian manifold be a fibre bundle\fg}, \emph{Proc. of Am. Math. Soc.}
  \textbf{11} (1960), p.~236--242.

\bibitem[Jam03]{ja03}
{\scshape P.~Jammes} -- {\og Sur le spectre des fibr\'es en tore qui
  s'effondrent\fg}, \emph{Manuscripta Math.} \textbf{110} (2003), no.~1,
  p.~13--31.

\bibitem[Jam05]{ja05}
\bysame , {\og Effondrements et petites valeurs propres des formes
  diff\'erentielles\fg}, \emph{S\'emin. th\'eor. spectr. g\'eom.} \textbf{23}
  (2005), p.~115--124.

\bibitem[Jam10]{ja10}
\bysame , {\og Effondrement, spectre et propriétés diophantiennes des flots
  riemanniens\fg}, \emph{Ann. inst. Fourier} \textbf{60} (2010), no.~1,
  p.~257--290.

\bibitem[Jam11]{jahdr}
\bysame , {\og Autour de la g\'eom\'etrie du laplacien de
  {H}odge-de~{R}ham\fg}, Habilitation à diriger des recherches, universit\'e
  Nice Sophia Antipolis, d\'ec. 2011.

\bibitem[Li80]{li80}
{\scshape P.~Li} -- {\og On the {S}obolev constant and the $p$-spectrum of a
  compact riemannian manifold\fg}, \emph{Ann. scient. \'Ec. norm. sup. (4)}
  \textbf{13} (1980), p.~451--469.

\bibitem[Lot02]{lo02}
{\scshape J.~Lott} -- {\og Collapsing and the differential form {L}aplacian:
  the case of a smooth limit space\fg}, \emph{Duke Math. Journal} \textbf{114}
  (2002), no.~2, p.~267--306.

\bibitem[LY80]{ly80}
{\scshape P.~Li {\normalfont \smfandname} S.~Yau} -- {\og Estimates of
  eigenvalues of a compact riemannian manifold\fg}, in \emph{Proceedings
  Symposium on Pure Math.}, vol.~36, 1980, p.~205--239.

\bibitem[Man08]{ma08}
{\scshape T.~Mantuano} -- {\og Discretization of {R}iemannian manifolds applied
  to the {H}odge {L}aplacian\fg}, \emph{Amer. J. Math.} \textbf{130} (2008),
  no.~6, p.~1477--1508.

\bibitem[Ram05]{ra05}
{\scshape S.~Ramanan} -- \emph{Global calculus}, Graduate Studies in
  Mathematics, vol.~65, Amer. Math. Soc., 2005.

\bibitem[Sch80]{sc80}
{\scshape W.~Schmidt} -- \emph{Diophatine approximations}, Lecture notes in
  mathematics, vol. 785, Springer Verlag, 1980.

\end{thebibliography}

\end{document}